\documentclass[11pt]{article}
\usepackage{a4wide,amsfonts,amssymb,amsmath}
\def\proof{{\bf Proof\quad}}
\def\beginpf{\proof}
\def\qed{\hfill\rule{2.2mm}{2.2mm}\vspace{1ex}}
\def\endpf{\qed}

\newtheorem{theorem}{Theorem}[section]

\newtheorem{corollary}[theorem]{Corollary}
\newtheorem{definition}[theorem]{Definition}

\newtheorem{lemma}[theorem]{Lemma}

\newtheorem{remark}[theorem]{Remark}

\def\trace{\mathop{\rm tr}\nolimits} 

\def\eps{\varepsilon}
\def\w{\omega}
\def\cZ{\mathcal Z}
\def\CC{\mathbb C}
\def\JJ{\mathbb J}

\def\DD{\mathbb D}

\def\NN{\mathbb N}
\def\RR{\mathbb R}
\def\TT{\mathbb T}
\def\LL{\mathcal L}
\def\cJ{\mathcal J}

\def\cR{\mathcal R}
\def\cU{\mathcal U}

\def\wB{\widetilde{\mathcal B}}
\def\eins{\mathbf 1}
\def\ds{\displaystyle}
\newcommand{\tr}{\mathop{\rm tr}\nolimits}

\newcommand{\HH}{\mathcal H}
\newcommand{\KK}{\mathcal K}

\newcommand{\carl}{\mathop{\rm Carl}\nolimits}

\newcommand{\range}{\mathop{\rm range}\nolimits}
\newcommand{\empset}{\emptyset }
\newcommand{\re}{\mathop{\rm Re}\nolimits}

\newcommand{\rank}{\mathop{\rm rank}\nolimits}
\renewcommand{\span}{\mathop{\rm span}\nolimits}

\def\text{\mbox}
\def\kasten{\endpf}

\title{Tangential interpolation in weighted vector-valued $H^p$ spaces, with applications}
\author{{Birgit Jacob\thanks{Department of Applied Mathematics, Delft University of Technology, P.O.~Box 5031, 2600 GA Delft, The Netherlands, \tt b.jacob@tudelft.nl}
\qquad Jonathan R.~Partington\thanks{School of Mathematics,
University of Leeds,
Leeds LS2 9JT, U.K. \tt
J.R.Partington@leeds.ac.uk}
\qquad Sandra Pott\thanks{Department of Mathematics, University of Glasgow,
Glasgow, G12 8QW, U.K. \tt s.pott@math.gla.ac.uk}}}

\begin{document}
\maketitle

\begin{abstract}
In this paper, norm estimates are obtained for the
problem of minimal-norm tangential interpolation by vector-valued analytic functions in weighted $H^p$ spaces,
 expressed
in terms of the Carleson constants of related scalar measures. Applications are given
to the notion of $p$-controllability properties of linear semigroup
systems and controllability by functions in certain Sobolev spaces.
\end{abstract}


{\bf 2000 Subject Classification.} 30D55, 30E05, 47A57, 47D06, 93B05.

\section{Introduction and Notation}
Given a Hilbert space $\HH$, operators $G_1,\ldots,G_n$
on $\HH$, vectors $a_1,\ldots, a_n$ in $\HH$, and
$z_1,\ldots,z_n$ in $\CC_+$ we estimate the minimal norm
of a function $f \in H^p_W(\CC_+,\HH)$, $1 \le p < \infty$, satisfying the
interpolation conditions
\[
G_k f(z_k)=a_k \quad (k=1,\ldots,n)
\]
(all necessary notation is explained below).
This can be regarded as a
problem of {\em tangential interpolation}
in the sense of \cite{BGR}.
We shall see that in many cases a sharp estimate can be given in terms
of the Carleson constants of various scalar measures.

In the paper \cite{jpp06}, certain weighted vector-valued
generalizations of the Shapiro--Shields interpolation theory \cite{ss} for the
Hardy space $H^2$ of the right-hand complex half-plane $\CC_+$ were
achieved. The central tool was a modification of an approach of
McPhail \cite{mcphail} to the matrix case via matrix
Blaschke--Potapov products (see e.g.~\cite{Pel03}), which allowed a
unified treatment of tangential interpolation results in the
literature as well as their extension to the general weighted case
(in the sense of matrix weights in the target space).

The purpose of the present paper is to extend this weighted
tangential interpolation theory to interpolation by functions in
vector-valued $H^p$ spaces, $1 \le p < \infty$, on the right half
plane. In certain cases, we can also deal with matrix-weighted
vector-valued $H^p$ spaces.
Note that matrix weights appear in two
different meanings here. First, we use \emph{``matrix weights in the
target space''}, in the sense that given a sequence of distinct
points $(z_k)$ in the right half plane $\CC_+$ and a sequence of $N
\times N$ matrices $(G_k)$, thought of as weights, with ranges $J_k
\subseteq \CC^N$, we try to find for each sequence in $(a_k)$ in
$\ell^p(J_k)$ an interpolating $\CC^N$ valued function $f$ in an
appropriate space with $G_k f(z_k)= a_k$ for all $k$. Weights in
this sense are
 very useful for questions of controllability in linear systems
with multidimensional input space governed by diagonal semigroups,
as discussed in \cite{jpp06}. The main interesting case here is the
\emph{``tangential"} case, i.e., where $\rank G_k=1$ for all $k \in
\NN$.

Second, the weights newly introduced in the present paper are matrix weights on the space of interpolating functions.
Namely, rather than interpolating by functions in $H^2(\mathbb C_+, \CC^N)$, we will seek to interpolate by functions
in the weighted vector-valued $H^p$-space
\begin{multline*}
   H^p_W(\CC_+,\CC^N) 
    = \{ f: \CC_+ \rightarrow \CC^N \text{ analytic: } \sup_{\eps >0}
   \int_{i \RR} \langle W(t)^{2/p} f(t + \eps),f(t + \eps) \rangle^{p/2} dt < \infty \},
\end{multline*}
where $W$ is a measurable function on $i\RR$ taking values a.e.~in the positive invertible $N \times N$ matrices.
We want to refer to $W$ as a \emph{``matrix-weight in the function space''}. The motivation in this case is
given by questions of controllability by functions in certain Sobolev spaces, which are new even in the scalar
case.

Again, the approach of McPhail modified to the matrix case will play an important role, together with
the theory of matrix $A_2$ weights.\\

In Section \ref{sec:2} we give norm estimates for
the minimum-norm interpolation problem. Applications to various notions of controllability
are contained in Section \ref{sec:3}.\\

We shall frequently use the following notation. Let $(z_k)_{k \in
\NN}$ be a Blaschke sequence of pairwise distinct elements
 in the right half plane $\CC_+ = \{ z \in \CC: \re z >0 \}$.
Let $b_k(z) = \frac{z-z_k}{z + \bar z_k}$ denote the Blaschke factor
for $z_k$. For $n \in \NN$, $1 \le k \le n$, let $B_n(z) =
\prod_{j=1}^n b_j(z)$, $ B_{n,k}(z) = \prod_{j=1, j \neq k}^n
b_j(z)$, $b_{n,k}= B_{n,k}(z_k)$, $b_{\infty,k}= \lim_{n \to \infty}
B_{n,k}(z_k)$.

Also $k_{z_k}= \frac{1}{2\pi} \frac{1}{z + \bar z_k}$ denotes the
reproducing kernel at $z_k$, so that $\langle
f,k_{z_k}\rangle=f(z_k)$ for all $f \in H^2(\CC_+)$.

For an index $p$ with $1 \le p < \infty$, we use $p'$ to denote the conjugate index $p/(p-1)$.

\section{Interpolation} \label{sec:2}
The aim of this section is to extend the results of McPhail
\cite{mcphail} to a vector setting. We begin by collecting some
tools.

\subsection{Carleson--Duren type
embedding theorems for matrix measures} \label{ssec:carl} 
Recall the
classical Carleson--Duren Embedding Theorem (see e.g.~\cite{Nik86}, \cite{duren}).
\begin{theorem} \label{thm:ccarleson}
Let $\mu$ be a non-negative Borel measure on the  right half plane
$\CC_+$ and let $1\le \alpha < \infty$. Then the following are equivalent:
\begin{enumerate}
\item The embedding
$$
     H^p(\CC_+) \rightarrow L^{\alpha p}( \CC_+, \mu)
$$
is bounded for some (or equivalently, for all) $1 \le p < \infty$.
\item
There exists a constant $C> 0$ such that
$$
    \int_{\DD} | k_\lambda(z)|^2 d \mu(z) \le C \|k_\lambda\|_{H^{2\alpha}}^2 \text{ for all } \lambda \in \CC_+
$$
\item
$$
    \mu(Q_I) \le C |I|^\alpha \text{ for all intervals } I \subset \RR,
$$
where $Q_I = \{ z= x + iy \in \CC_+: y \in I, 0 < x < |I|\}$.
\end{enumerate}
In this case, $\mu$ is called a $\alpha$-Carleson measure.
\end{theorem}

For $0 < \alpha < 1$, we will call $\mu$ an $\alpha$-Carleson measure if the embedding
$$
     H^p(\CC_+) \rightarrow L^{\alpha p}( \CC_+, \mu)
$$
is bounded for some (or equivalently, for all) $1 \le p < \infty$.  In this case, conditions (2) and (3) of the Theorem
are no longer sufficient to $\mu$ to be $\alpha$-Carleson, but they are easily seen to be necessary.
A necessary and sufficient condition for the case $ \alpha <1$ can be found in \cite[Thm.~C]{Lue91},  and will be summarized in the following theorem.

For a scalar or operator valued regular Borel measure $\mu$ on $\CC_+$, let $S_\mu$ denote the \emph{balayage}
of $\mu$,
$$
    S_\mu(i \w) = \int_{\CC_+} p_z(i \w) d\mu(z),
$$
where  
\begin{equation}\label{poisson}
p_z(i\omega)={\pi}^{-1}\frac{\ds x}{\ds x^2+(y-\omega)^2}
\end{equation}
 denotes the Poisson kernel for $z=x+iy$ on $i \RR$.

\begin{theorem} \cite{Lue91}
Let $0< \alpha <1$ and let $\mu$ be a (scalar-valued) non-negative
regular Borel measure on $\CC_+$.
Then   
$$
H^p(\CC_+) \rightarrow L^{\alpha p}( \CC_+, \mu)
$$
is bounded for some, and equivalently, for all $0 < p < \infty$,
if and only if $S_\mu \in L^{1/(1-\alpha)}(i \RR)$.
\end{theorem}

For a discrete measure with finite support,
 $\mu =\sum_{k=1}^N A_k \delta_{z_k}$ on $\CC_+$, the balayage (if it exists)
can be conveniently expressed as
\begin{equation*}
 S_\mu(i \w) = \sum_{k=1}^N A_k p_{z_k}(i \w).
\end{equation*}

We can trivially include the notion of a $0$-Carleson measure here,
denoting a finite measure, and find that
\begin{equation} \label{eq:zerocarl}
H^\infty(\CC_+) \rightarrow L^{p}( \CC_+, \mu)
\end{equation}
is bounded for some, and equivalently, for all $0 < p < \infty$,
if and only if $\mu$ is $0$-Carleson.

We write $\carl_\alpha(\mu)$ for the infimum of constants satisfying \ref{thm:ccarleson} (2) in case
$\alpha \ge 1$, respectively $\|S_\mu\|_{1/(1- \alpha)}$ in case $0< \alpha <1$. With this notation,
the known results yield easily that
$$
\|H^p(\CC_+) \rightarrow L^{\alpha p}( \CC_+, \mu)\| \approx \carl_\alpha^{1/(\alpha p)}(\mu)
$$
for $1\le p < \infty$, $\alpha >0$, with equivalence constants depending only on $p$ and $\alpha$.

Let us use the following notation: For $0 < p < \infty$,
$\HH$ a finite or infinite-dimensional Hilbert space,
$$
    H^p(\CC_+, \HH) = \left\{f: \CC_+ \rightarrow \HH \text{ analytic: }
     \sup_{\eps >0} \int_{i \RR}\|f(t + \eps)\|^p \, d t < \infty   \right\},
$$
and
$$
    H^p(\CC_-, \HH) = \left\{f: \CC_+ \rightarrow \HH \text{ anti-analytic: }
     \sup_{\eps >0} \int_{i \RR}\|f(t + \eps)\|^p \, d t < \infty   \right\}.
$$
Although a full operator analogue of even the classical Carleson
Embedding Theorem is not known, the following
 is easily proved.

\begin{theorem}
\label{thm:opcarleson}
 Let $\mu$ be a non-negative operator-valued Borel measure on the  right half plane $\CC_+$.
For $0 < p < \infty$, let
$$
L^p(\CC_+, \HH, \mu) = \left\{ f: \CC_+ \rightarrow \HH \text{ strongly
measurable}: \int_{\DD}  \langle d\mu(z)^{2/p} f(z), f(z)
\rangle^{p/2} < \infty\right\}.
$$

Let $\|\mu\|$ be the total variation of $\mu$,
$$
   \| \mu\|(A) = \sup \left\{ \sum_{i=1}^n \|\mu(A_i)\| :
     A_1, \dots, A_n \text{ pairwise disjoint, } A_1 \cup \dots \cup A_n = A, n \in \NN \right\}.
$$
Suppose that $\|\mu\|$ is a scalar $\alpha$-Carleson measure.

   Then
the embedding
$$
     H^p(\CC_+,\HH ) \rightarrow L^{\alpha p}( \CC_+,\HH, \mu)
$$
is bounded for $1\le p < \infty$, $0< \alpha < \infty$, and
the embedding
$$
     H^\infty(\CC_+,\HH ) \rightarrow L^{p}( \CC_+,\HH, \mu)
$$
is bounded for $0< p < \infty$, $\alpha =0$.

If $\dim \HH < \infty$, then the reverse is also true.
\end{theorem}
\proof A proof is stated here for the convenience of the reader.
Let $f \in H^p(\CC_+, \HH)$.
Choosing an orthonormal basis $(e_j)$ of $\HH$ and writing $f_j = \langle f(\cdot), e_j \rangle$, we obtain
\begin{eqnarray*}
    \int_{\CC_+} \langle d \mu(z)^{2/(\alpha p)} f(z), f(z) \rangle^{\alpha p/2} &\le &
    \int_{\CC_+} (\sum_{j=1}^\infty  |f_j(z)|^2)^{\alpha p/2} d \|\mu\|(z)  \\
   & \approx& \int_{\CC_+} \int_0^1 |\sum_{j=1}^\infty  r_j(s) f_j(z)|^{\alpha p} ds  d \|\mu\|(z)  \\
   &\lesssim &  \left(\int_{0}^1 \int_{i\RR} |\sum_{j=1}^\infty r_j(s) f_j(t)|^p dt ds\right)^\alpha \\
   & \approx & \left(\int_{i \RR} (\sum_{j=1}^\infty |f_j(t)|^2)^{p/2} dt\right)^\alpha \\
   &=& \|f\|^{\alpha p}_{H^p}
\end{eqnarray*}
by the Carleson--Duren Theorem, respectively the definition of
a $\alpha$-Carleson measure, in the scalar case. Here, the $r_j$, $j \in \NN$, denote the Rademacher functions on
$[0,1]$, and we use Khintchine's inequalities in lines 2 and 4, and constants depend only on $p$.

For the reverse implication in the finite-dimensional case, just
note that a comparison of trace and operator norm gives that
$\|\mu\|$ is a scalar $\alpha$-Carleson measure if and only if $\trace \mu$
is $\alpha$-Carleson.  Let $e_1, \dots, e_N$ denote an orthonormal basis of
$\HH$. Then in case $0< \alpha < \infty$, the reverse implication follows easily from the identity
$$
    \int_{\CC_+} |f(z)|^{\alpha p} d \trace \mu = \sum_{i=1}^N \int_{\CC_+} \langle d
    \mu(z)^{2/(\alpha p)} f(z) e_i, f(z) e_i \rangle^{\alpha p/2} \quad (f \in
    H^p(\CC_+) )
$$
and the scalar case. In case $\alpha =0$, apply boundedness of the
embedding  $H^\infty(\CC_+,\HH ) \rightarrow L^{p}( \CC_+,\HH, \mu)$
to the constant $\HH$-valued functions $e_1, \dots, e_N$.
\qed

Finally, we want to deal with a less trivial case, the case of matrix-weighted embeddings.

A strongly measurable function $ W: i\RR \rightarrow \LL(\HH)$ is called an \emph{operator weight}, if it takes
values in the positive invertible operators in $\LL(\HH)$
a.e. and $\|W\|, \|W^{-1}\| \in L^1_{loc}(i\RR)$.
If $\dim \HH < \infty$, we speak of matrix weights.
A matrix weight $W$  is called matrix $A_2$, if
\begin{equation} \label{eq:a2}
      \sup_{I \subset \mathbb R, I \text{ bounded interval}}
     \left\| \left(\frac{1}{|I|}\int_I W^{-1}(t)\,dt\right)^{1/2}\left(\frac{1}{|I|}\int_I W(t)\,dt\right)\left(\frac{1}{|I|}\int_I W^{-1}(t)\,dt\right)^{1/2}\right\| < \infty
\end{equation}
(see \cite{tv1}). An equivalent formulation is the ``invariant
matrix $A_2$ condition''
\begin{equation}   \label{eq:inva2}
 \sup_{z \in \CC_+}
     \| (W^{-1}(z))^{1/2} W(z)  (W^{-1}(z))^{1/2}\| < \infty,
 \end{equation}
(see \cite[Lem.~2.2]{tv1}), where $W(z)$, $W^{-1}(z)$ denote the
Poisson extension of $W$ resp. $W^{-1}$ in $z \in \CC_+$ (so in
general $W^{-1}(z) \neq (W(z))^{-1}$). The invariant $A_2$ condition implies at once
\begin{equation*}
         W^{-1}(\cdot) \approx (W(\cdot))^{-1} .
\end{equation*}

For matrix $A_2$ weights, we have the following matrix-weighted embedding theorem, essentially taken
from \cite{tv2}:
\begin{theorem}
\label{thm:weightcarleson}
 Let $\mu$ be a non-negative $N \times N$ matrix-valued Borel measure on the  right half plane $\CC_+$,
and let $W$ be a matrix $A_2$ weight on $i\RR$. We write $W(z)$ for the harmonic extension of $W$
in $z \in \CC_+$. Let
$$
L^2_W(i \RR, \CC^N) = \left\{ f: \CC_+ \rightarrow \CC^N \text{ measurable}:
\int_{i \RR}  \langle W(t) f(t), f(t) \rangle dt< \infty\right\},
$$
\begin{multline*}
L^2_W(\CC_+, \CC^N, d \mu) = \\
\left\{ f: \CC_+ \rightarrow \CC^N \text{ measurable}:
\int_{\CC_+}  \langle(W(z))^{1/2} d\mu(z) (W(z))^{1/2} f(z), f(z) \rangle < \infty\right\}.
\end{multline*}

Then
the embedding
$$
     L^2_W(i \RR,\CC^N ) \rightarrow L^2_W( \CC_+,\CC^N, d\mu)
$$
is bounded, if and and only if $\mu$ is a matrix Carleson measure.
\end{theorem}
\proof ``$\Leftarrow$'' For the case of the unit disk and scalar
measure $\mu$, this is proved in \cite[Lem.~4.1]{tv2}. The case of
the right half plane for scalar measures is proved similarly. To
obtain the boundedness of the embedding for a matrix measure $\mu$,
just note that $\tr \mu$ is a Carleson measure by Theorem
\ref{thm:opcarleson} and that
$$
\int_{\CC_+}  \langle(W(z))^{1/2} d\mu(z) (W(z))^{1/2} f(z), f(z) \rangle
\le
\int_{\CC_+}  \langle W(z)f(z), f(z) \rangle d \tr \mu(z).
$$

``$\Rightarrow$'' As in the case of the unit disk, the matrix $A_2$ condition implies a certain factorization
of the weight.
Namely, there exist matrix-valued functions $F, G \in L^2_{loc}(i \RR, \CC^{N \times N})$ such that
$F \circ M$, $G \circ M$ are outer functions in $H^2( \DD, \CC^{N \times N})$, and
$W = F^* F$, $W^{-1} = G G^*$  a.e. on $i \RR$. Here, $M: \DD \rightarrow \CC_+$ denotes the Cayley transform
$ z \mapsto \frac{1-z}{1 +z}$.
(For the case of the unit disk, see the proof of \cite[Thm.~3.2]{tv2}. The case of the right half-plane $\CC_+$
then follows easily from the fact that composition with $M$ maps a matrix $A_2$ weight on $i \RR$
to a matrix $A_2$ weight on $\TT$.)

As in the proof of \cite[Thm.~3.2]{tv2},
$$
    | \det G(z) | = |\det F(z)|^{-1}   \quad \text{ for } z \in \CC_+
$$
and
\begin{equation}   \label{eq:csdir}
   F^*(z) F(z) \le (F^* F)(z), \quad
G(z) G^*(z) \le (G G^*)(z) \text{ for all } z \in \CC_+
\end{equation}
by Cauchy--Schwarz. Since by the matrix $A_2$ condition there exists a constant $C >0$ with
$$
   (W^{-1}(z))^{1/2} W(z)  (W^{-1}(z))^{1/2}  \le C^2 \eins \text{ for all } z \in \CC_+,
$$
it follows that
$$
    \sup_{z \in \CC_+} |\det(W(z))| |\det(W^{-1}(z))| =
        \sup_{z \in \CC_+} \frac{|\det(W(z))|}{|\det F(z)|^2} \frac{|\det(W^{-1}(z))|}{|\det G(z)|^2} \le C^{2N},
$$
with both factors bounded by below $1$ because of (\ref{eq:csdir}). Thus
$$
     |\det G(z)| \le (\det W^{-1}(z))^{1/2} \le C^N |\det G(z)| \quad \text{ for all } z \in \CC_+.
$$
As
$$
    0 \le (W^{-1}(z))^{-1/2} G(z) G^*(z) (W^{-1}(z))^{-1/2} \le \eins
$$
and
$$
\frac{1}{C^{2N}} \le \det((W^{-1}(z))^{-1/2} G(z) G^*(z) (W^{-1}(z))^{-1/2}) \le 1 \text{ for } z \in \CC_+,
$$
it follows that
\begin{equation}     \label{eq:winvest}
 \frac{1}{C^{2N}} \eins \le (W^{-1}(z))^{-1/2} G(z) G^*(z) (W^{-1}(z))^{-1/2} \le \eins \text{ for } z \in \CC_+.
\end{equation}
Applying the invariant matrix $A_2$ condition yet again,
\begin{equation}     \label{eq:wwinvest}
 \frac{1}{C^{2N}} \eins \le (W(z))^{1/2} G(z) G^*(z) (W(z))^{1/2} \le C^2 \eins \text{ for } z \in \CC_+.
\end{equation}
Now let
$\lambda \in \CC_+$, and let $e_1, \dots, e_N$  be the standard basis
$\CC^N$. Then (\ref{eq:wwinvest}) implies
\begin{eqnarray*}
   && \frac{1}{C^{2N}} \int_{\CC_+}  |k_\lambda(z)|^2 d \tr \mu(z) \\
   &\le& \int_{\CC_+}  \tr ( G(z)^* (W(z))^{1/2} d \mu(z) (W(z))^{1/2}G(z)) |k_\lambda(z)|^2 \\
   &=&  \sum_{i=1}^N \int_{\CC_+}
    \langle  (W(z))^{1/2} d \mu(z) (W(z))^{1/2}G(z) k_\lambda(z)e_i, G(z)k_\lambda(z)e_i \rangle \\
  &\le&
    {\tilde C}^2 \int_{i \RR} \sum_{i=1}^N \langle W(t) G(t) k_\lambda e_i, G(t) k_\lambda e_i \rangle dt \\
   &=& {\tilde C}^2 \int_{i \RR}  \tr (G^*(t) W(t) G(t))  |k_\lambda(t)|^2 dt
    = N {\tilde C}^2 \|k_\lambda\|^2_{H^2},
\end{eqnarray*}
where $\tilde C$ denotes the norm of the embedding $L^2_W(i
\RR,\CC^N ) \rightarrow L^2_W( \CC_+,\CC^N, d\mu)$. Thus $\tr \mu$
is a Carleson measure, and $\mu$ is a matrix Carleson measure by
Theorem \ref{thm:opcarleson}. \qed

\subsection{Certain shift-invariant subspaces in $H^p(\CC_+, \HH)$}

Let $\HH$ be a separable Hilbert space.
 Exactly as in Lemma 2.4 in
\cite{jpp06}, one proves
\begin{lemma} \label{lemm:blaschke} Let $1 \le p < \infty$, and
let $(z_k)_{k \in \NN}$ be a sequence of pairwise distinct
elements of $\CC_+$. For each $k \in \NN$, let
$L_k \subseteq \HH$ be a closed linear subspace of $\HH$.
Let
$\LL_n = \{ f \in H^p(\CC_+, \HH): f(z_k) \in L_k \text{ for }1 \le k \le n\}$.
Then $\LL_n=\Theta^L_n H^p(\CC_+, \HH)$, where $\Theta^L_n$ is the matrix-valued
Blaschke--Potapov product
\begin{equation*}
  \Theta_n^L(z) = (b_1(z)P^{\perp}_{\tilde L_1} + P_{\tilde L_1})\cdots
                                         (b_n(z)P_{\tilde L_n}^\perp + P_{\tilde L_n}) \quad (z \in\CC_+),
\end{equation*}
where
$$
\tilde L_1 = L_1, \quad \tilde L_k = \Theta^L_{k-1}(z_k)^{-1} L_k \quad\text{ for } 2 \le k \le n
$$
and $P_{\tilde L_k}$ is the orthogonal
projection $\HH \rightarrow \tilde L_k$.

One sees easily that $P_{L^\perp_k}\Theta^L_n(z_k) =0$ for $k=1, \dots, n$.

If $(z_k)$ is a Blaschke sequence, then $\Theta^L_n$ converges
normally (uniformly on compact subsets of $\CC_+$) to an inner
function $\Theta^L$ with $\Theta^L H^p(\CC_+, \HH) =\cap_{n \in \NN}
\LL_n$.
\end{lemma}


\subsection{Interpolation Theorems}
With the notation of the $\Theta^L_{n}$, we can formulate our generalizations of McPhail's result
\cite[Thm.~2~(B)]{mcphail}. Let $W$ be an operator weight such that there exists $M \in \NN$ with
$$
     \int_{i \RR} \frac{1}{(1 + t)^M} (\|W(t)\|+ \|W^{-1}(t)\|) dt < \infty.
$$
We define for $1 \le p < \infty$
$$
  H^p_W(\CC_+, \HH) =
     \left\{ f: \CC_+ \rightarrow \HH \text{ analytic: } \sup_{\eps >0}
   \int_{i \RR} \langle W(t)^{2/p} f(t + \eps),f(t + \eps) \rangle^{p/2} dt < \infty \right\},
$$
\begin{multline*}
H^p_W(\CC_-, \HH) = \\
     \left\{ f: \CC_+ \rightarrow \HH \text{ anti-analytic: } \sup_{\eps >0}
   \int_{i \RR} \langle W(t)^{2/p} f(t + \eps),f(t + \eps) \rangle^{p/2} dt < \infty \right\}
\end{multline*}
and for $p= \infty$,
$$
  H^\infty_W(\CC_+, \HH) =
     \{ f: \CC_+ \rightarrow \HH \text{ analytic: } \sup_{\eps >0} \sup_{t \in i \RR}
       \| W(t) f(t + \eps )\| < \infty \},
$$
$$
H^\infty_W(\CC_-, \HH) =
     \{ f: \CC_+ \rightarrow \HH \text{ anti-analytic: } \sup_{\eps >0}
         \sup_{t \in i \RR} \| W(t) f(t + \eps )\|  < \infty \}.
$$

For $N= \dim \HH < \infty$, it was shown in \cite{tv1} that $W$ is a
matrix $A_2$ weight as in (\ref{eq:a2})  if and only if
\begin{equation}   \label{eq:tv}
       L^2_W(i \RR, \HH) \simeq H^2_W(\CC_+, \HH) \oplus H^2_W(\CC_-,
       \HH),
\end{equation}
with equivalence constants of norms only depending on $N$ and the
$A_2$ constant of $W$.

 With the above notation, we have the following duality relations for $1 < p <
 \infty$:
\[
 (L_W^p(i \RR, \HH) / H_W^p(\CC_+, \HH))^* = H_{W^{-1/{(p-1)}}}^{p'}(\CC_-, \HH)  =
\left\{ \bar f : f \in H_{\overline{W^{-1/{(p-1)}}}}^{p'}(\CC_+, \HH)\right\},\]
 where $\bar f$ stands for the coordinatewise
complex conjugate with respect to some fixed orthonormal basis of $\HH$, and $\overline{W^{-1}}$ stands
for the entry-wise complex conjugate of the matrix representation of $W^{-1}$ with respect to the chosen basis.
The duality is given by
\begin{multline*}
\qquad\qquad   \langle [f], g \rangle = \int_{i \RR} \langle f(t), g(t) \rangle_\HH dt \\
 \text{ for } [f] \in
       L_W^p(i \RR, \HH) / H_W^p(\CC_+, \HH) \,\,\mbox{ and}\,\, g \in H_{W^{-1/{(p-1)}}}^{p'}(\CC_-, \HH).
\end{multline*}
for $1 < p < \infty$. For $p=1$ we have
$ (L_W^1(i \RR, \HH) / H_W^1(\CC_+, \HH))^* = H^\infty_{W^{-1}}(\CC_{-}, \HH)$. \\

Let $\HH$ be a separable Hilbert space, and let $(G_k)_{k \in \NN}$
be a sequence of non-zero bounded linear operators on $\HH$ with
closed range. We will be particularly interested in the case of
finite-dimensional $\HH$ and of $G_k$ being of finite rank,
specifically of rank $1$. We write $I_k = (\ker G_k)^\perp$, $J_k =
\range G_k$ for $k \in \NN$, and denote $\dim I_k = \dim J_k$ by
$d_k$ in the case that $G_k$ has finite rank. In the finite rank
case, we write, slightly
 abusing notation, $G_k^{-1}: \CC^N \rightarrow I_k \subseteq \CC^N$
for the linear operator defined by $(G_k|_{I_k \to
J_k})^{-1}P_{J_k}$.
 We fix the Blaschke
sequence $(z_k)_{k \in \NN}$
 of pairwise distinct elements of $\CC_+$ and the sequence $(G_k)_{k \in \NN}$.

For $n \in \NN$, $1 \le p < \infty$, and $1\le s<\infty$, let
\begin{multline*}
   m_{n,p,s,W} =
\sup_{ a \in \bigoplus_{k =1}^n J_k, \|a\|_s \le 1}
        \inf \{ \|f\|_p: f \in H_W^p(\CC_+, \HH),\,\, G_k f(z_k) = a_k,\,\, k=1,\dots, n \}.
\end{multline*}
and
\begin{multline*}
   m_{p,s,W} =
\sup_{ a \in \bigoplus_{k =1}^\infty J_k, \|a\|_s \le 1}
        \inf \{ \|f \|_p: f \in H_W^p(\CC_+, \HH),\,\, G_k f(z_k) = a_k,\,\, k\in \NN \}.
\end{multline*}

A weak$^*$ compactness argument shows that $m_{p,s,W} = \sup_{n \in
\NN} m_{n,p,s,W}$.

Here comes the main interpolation result.

\begin{theorem} \label{thm:carlesone} Let $1 \le p < \infty$.
Let $(G_k)_{k \in \NN}$, $(I_k)_{k \in \NN}$, $(z_k)_{k \in \NN}$,
$m_{p,s,W}$ be defined as above. Let $\Theta^{I^\perp}$
 be the inner function associated to the sequence
$(z_k)_{k \in \NN}$
and the sequence of subspaces $(I_k^\perp)_{k \in \NN}$ as in Lemma \ref{lemm:blaschke}.
\begin{enumerate}
\item If $1\le p < \infty$ and $1 < s < \infty$, then
$$
 m_{p,s,W} = \|  \cJ\|_{H^{p'}_{\tilde W_{p}}(\CC_-, \HH) \rightarrow L^{s'}(\CC_+, \HH, d \mu_{s})}
$$
where
$$
   \mu_{s} = \sum_{k=1}^\infty
     \frac{|2\re  z_k|^{s'}}{|b_{\infty,k}|^{s'}}  (\overline{\Theta^I(z_k)}^* G_k^{-1}(G_k^{-1})^*
                                                          \overline{\Theta^I(z_k)})^{s'/2} \delta_{z_k},
$$
$$
   \tilde W_{p} = {\Theta^{I^\perp}}^* W^{-1/(p-1)}\Theta^{I^\perp}
   \text{ if } 1< p < \infty, \; \tilde W_{1} = {{\Theta^{I^\perp}}^* W^{-1}\Theta^{I^\perp}},
$$
and $\cJ$ is the natural embedding operator.

\item If $1 \le p < \infty$, $s=1$, then
$$
 m_{p,1,W} = \|  \cJ_{{1}}\|_{H^{p'}_{\tilde W_{p}}(\CC_-, \HH) \to \ell^\infty(\HH)} \quad(n \in \NN),
$$
where $\tilde W_p$ is as above,
and $\cJ_{1}$ is the embedding
$$
  \cJ_{1}:H_{\tilde W_{p}}^{p'}(\CC_-, \HH) \rightarrow \ell^\infty(\HH),
                       \quad f \mapsto (2 \frac{\re z_k}{|b_{\infty,k}|} (G_k^{-1})^*\overline{\Theta^{I}(z_k)} f(z_k) )_{k \in \NN}.
$$
\end{enumerate}
\end{theorem}

\proof As in the case $p=2$ \cite{jpp06}, we prove this by first interpolating finitely many points
 and then using the uniform convergence of the Blaschke products
$\Theta_n^{I}$, $\Theta_n^{I^\perp}$ on compact subsets of $\CC_+$:
\begin{lemma} \label{lemma:carlfinite} Let $1 \le s, p < \infty$.
Let $(G_k)_{k \in \NN}$, $(I_k)_{k \in \NN}$, $(z_k)_{k \in \NN}$,
$(m_{n,p,s,W})_{n \in \NN}$ be defined as above. Let $\Theta_n^{I^\perp}$
 be the inner functions associated to the tuple
$z_1, \dots, z_n$
and the subspaces $I_1^\perp, \dots, I_n^\perp$ as in Lemma \ref{lemm:blaschke}.
\begin{enumerate}
\item If $1<s< \infty$, $1 \le p < \infty$, then
$$
 m_{n,p,s,W} = \|  \cJ \|_{H^{p'}_{\tilde W_{n,p}}(\CC_-, \HH) 
                                          \rightarrow L^{s'}(\CC_+, \HH, d \mu_{n,s})}  \quad(n \in \NN),
$$
where
$$
   \mu_{n,s} = \sum_{k=1}^n
     \frac{|2\re  z_k|^{s'}}{|b_{n,k}|^{s'}}  (\overline{\Theta^I_n(z_k)}^* G_k^{-1}(G_k^{-1})^*
                                                          \overline{\Theta^I_n(z_k)})^{s'/2} \delta_{z_k},
$$
$$
   \tilde W_{n,p} = {\Theta_n^{I^\perp}}^* W^{-1/(p-1)}\Theta_n^{I^\perp}
   \text{ if } 1< p < \infty, \; \tilde W_{n,1} = {{\Theta_n^{I^\perp}}^* W^{-1}\Theta_n^{I^\perp}},
$$
and $\cJ$ is the natural embedding operator.

\item If $s=1$, $1 \le p < \infty$, then
$$
 m_{n,p,1,W} = \|  \cJ_{{n,1}}\|_{    H^{p'}_{\tilde W_{n,p}}(\CC_-,
 \HH)  \to \ell_n^\infty(\HH)}
  \quad(n \in \NN),
$$
where
$\tilde W_{n,p}$ is as above and
$\cJ_{{n,1}}$ is the embedding
$$
   \cJ_{n,1}:H_{\tilde W_{n,p}}^{p'}(\CC_-, \HH) \rightarrow \ell_n^\infty(\HH),
                       \quad f \mapsto (2 \frac{\re z_k}{|b_{n,k}|} G_k^{-1}\Theta^{I^\perp}_n(z_k) f(z_k) )_{k=1, \dots, n}.
$$
\end{enumerate}
\end{lemma}

\proof
 Choose $M \in \NN$
such that
$$
   \int_{i \RR} \frac{1}{(1 + t)^M} (\|W(t)\|+ \|W^{-1}(t)\|) dt < \infty.
$$

For $a \in \oplus_{k=1}^n J_k$, let
$$
    \Phi_a(z) = \sum_{k=1}^n b_k(z)^{-1} \frac{(1 + z_k)^M}{(1 + z)^M}  (\Theta^{I^\perp}_{n,k})^{-1} G_k^{-1} a_k \quad
(z\in \CC_+ \backslash\{z_1, \dots, z_n\}).
$$
(recall that $ G_k^{-1}: J_k \rightarrow I_k$, $(\Theta^{I^\perp}_{n,k})^{-1}: I_k \rightarrow \HH$).
By choice of $M$, $ \Phi_a \in L^1_W(i \RR, \LL(\HH)) \cap L^\infty(i \RR, \LL(\HH))$.
Let $F_a =  \Theta_n^{I^\perp} \Phi_a$.\\

As in Lemma 2.7 in \cite{jpp06}, one proves that $F_a$ extends to an
analytic function on $\CC_+$ and that $G_k F_a(z_k) = a_k$ for $k=1,
\dots, n$.

 Since $\Phi_a|_{i \RR} \in L_W^1(i \RR, \LL(\HH))
\cap L^\infty(i \RR, \LL(\HH)) $ and $\Theta_n^{I^\perp}$ is inner,
$F_a \in H_W^p(\CC_+, \CC^N)$. So $F_a$ is an interpolating function
 in the desired sense. We now seek to solve the minimal-norm
interpolation problem.

 For all $g \in
H^p(\CC_+, \HH)$ with $G_k g(z_k)= a_k$, we have $g(z_k) - F_a(z_k)
\in I_k^\perp$.

By Lemma \ref{lemm:blaschke},
\begin{eqnarray}
   m_{n,p,s,W} &=& \sup_{a \in \oplus_{k=1}^n J_k, \|a\|_s=1}
      \inf_{f \in H^p_W(\CC_+, \HH)} \| F_a - \Theta_n^{I^\perp}f\|_{p,W} \nonumber\\
&=&  \sup_{a \in \oplus_{k=1}^n J_k, \|a\|_s=1} \,
      \inf_{f \in H^p_{{\Theta_n^{I^\perp}}^* W\Theta_n^{I^\perp}}(\CC_+, \HH)}
         \| \Phi_a - f\|_{p,{\Theta_n^{I^\perp}}^* W\Theta_n^{I^\perp}} \nonumber \\
& =&
        \sup_{a \in \oplus_{k=1}^n J_k, \|a\|_s=1}
       \| [\Phi_a]\|_{ L^p_{{\Theta_n^{I^\perp}}^* W\Theta_n^{I^\perp}}(i \RR, \HH) /
     H^p_{{\Theta_n^{I^\perp}}^* W\Theta_n^{I^\perp}}(\CC_+, \HH)}
      \nonumber \\
&=&  \sup_{a \in \oplus_{k=1}^n J_k, \|a\|_s=1}
    \sup_{ f \in H^{p'}_{\tilde W_{n,p}}(\CC_-, \HH), \|f\|=1}
       |\langle \Phi_a,f \rangle|, \nonumber \\
\end{eqnarray}
where $ \tilde W_{n,p} = {{\Theta_n^{I^\perp}}^*
  W^{-1/(p-1)}\Theta_n^{I^\perp}}$ for $1< p < \infty$ and $ \tilde W_{n,1} = {{\Theta_n^{I^\perp}}^* W^{-1}\Theta_n^{I^\perp}}$.
Now we have to distinguish between the cases $1 < s < \infty$ and $s=1$.  

 Temporarily writing $\cZ$ for \newline $H^{p'}_{\tilde
  W_{n,p}}(\CC_-, \HH) = ( L^p_{{\Theta_n^{I^\perp}}^* W\Theta_n^{I^\perp}}(i \RR, \HH) /
     H^p_{{\Theta_n^{I^\perp}}^* W\Theta_n^{I^\perp}}(\CC_+, \HH))^*
  $, we have for
$1< s < \infty$:
\begin{eqnarray*}
m_{n,p,s,W}&=&
        \sup_{a \in \oplus_{k=1}^n J_k, \|a\|_s=1}
    \sup_{ f \in \cZ, \|f\|=1}
       |\langle \Phi_a,f \rangle| \\
&=& \sup_{ f \in \cZ, \|f\|=1}
 \:  \sup_{a \in \oplus_{k=1}^n J_k, \|a\|_s=1}
              |\langle \Phi_a,  f\rangle| \\
&=& \sup_{f \in \cZ, \|f\|=1} \,  
     \sup_{a \in \oplus_{k=1}^n J_k, \|a\|_s=1}
 |\langle \sum_{k=1}^n  \overline{b_k(z)}   \left(\frac{1+z_k}{1+z}\right)^M  (\Theta_{n,k}^{I^\perp})^{-1} G_k^{-1} a_k,f\rangle|
\\
&=& \sup_{f \in \cZ, \|f\|=1}  
  \sup_{a \in \oplus_{k=1}^n J_k, \|a\|_s=1} \nonumber \\
&& 2 \pi \left|\left \langle \sum_{k=1}^n  -(z+\bar z_k)  \left(\frac{1+z_k}{1+z}\right)^M 
                  \langle (\Theta_{n,k}^{I^\perp})^{-1} G_k^{-1} a_k, f(z)\rangle_{\HH}, 
     \frac{1}{2\pi} \frac{1}{z+ \bar z_k} \right\rangle_{H^2(\CC_+)}\right|
\end{eqnarray*}
Thus, 
\begin{eqnarray}
 m_{n,p, s,W}
&=& \sup_{f \in \cZ, \|f\|=1 }  \, \sup_{a \in \oplus_{k=1}^n J_k, \|a\|_s=1}
2 \pi |\langle \sum_{k=1}^n  2 \re (z_k)  
                  \langle f(z_k), (\Theta_{n,k}^{I^\perp})^{-1} G_k^{-1} a_k\rangle_{\HH}|.
\nonumber\\
&=& {2 \pi} \sup_{ f \in \cZ, \|f\|=1} \, 
                     \sup_{a \in \oplus_{k=1}^n J_k, \|a\|_s=1}
              |\sum_{k=1}^n \langle 2 \re (z_k) a_k, (G_k^{-1})^* ((\Theta^{I^\perp}_{n,k})^{-1})^*  f(z_k)\rangle| 
\nonumber\\
&=& {2 \pi}\sup_{f \in \cZ, \|f\|=1}
  \left(\sum_{k=1}^n \|  \frac{2 \re (z_k)}{b_{n,k}}
(G_k^{-1})^*  \overline{\Theta_n^I(z_k)}U_n f(z_k)\|^{s'}\right)^{1/s'} 
\nonumber\\
&=& {2 \pi} \sup_{f \in \cZ, \|f\|=1}
                          \left(\sum_{k=1}^n \| \frac{2 \re (z_k)}{b_{n,k}}  (G_k^{-1})^*
                                                           \overline{\Theta_n^I(z_k)} f(z_k)\|^{s'}\right)^{1/s'} 
\nonumber\\
&=& \| \cJ|_{ H^{p'}_{\tilde W_{n,p}}(\CC_-, \HH) \rightarrow L^{s'}(\CC_+, \HH, \mu_{n,s})}\|.
\label{eq:main} 
\end{eqnarray}
Here, $U_n$ is a suitably chosen unitary operator, using Lemma 2.6 from \cite{jpp06} in the third line of the proof.

If $s=1$, then
\begin{eqnarray*}
m_{n,p,1}&=&  \sup_{a \in \oplus_{k=1}^n J_k, \|a\|_1=1}
      \sup_{f \in H^{p'}_{\tilde W_{n,p}}(\CC_-, \HH), \|f\|=1} |\langle \Phi_a,  f\rangle| \\
&=& \sup_{f \in  H^{p'}_{\tilde W_{n,p}}(\CC_-, \HH)              , \|f\|=1}
          \sup_{k \le n} \| 2 \re (z_k)  (G_k^{-1})^* ({\Theta^{I^\perp}_{n,k}}^{-1})^* f(z_k)\|
                                                                                                      \\
&=& \sup_{f \in H^{p'}_{\tilde W_{n,p}}(\CC_-, \HH)   , \|f\|=1}
          \sup_{k \le n} \| 2 \frac{\re (z_k)}{|b_{n,k}|}  (G_k^{-1})^*
      \overline{\Theta^I_n(z_k)} U_n {f}(z_k)\|\\
&=& \sup_{f \in H^{p'}_{\tilde W_{n,p}}(\CC_-, \HH) , \|f\|=1}
          \sup_{k \le n} \| 2 \frac{\re (z_k)}{|b_{n,k}|}  (G_k^{-1})^* \overline{\Theta^I_n(z_k)}  {f}(z_k)\|.
 \end{eqnarray*}

\qed

This finishes the proof of Theorem \ref{thm:carlesone}.
 \qed

In the unweighted case, we can instead consider a Carleson embedding for a simpler measure, restricted to
an invariant subspace of the shift operator:
\begin{corollary} \label{cor:inv} Let $1 \le  p < \infty$, $1< s < \infty$.
$$
 m_{n,p,s}
          =  
\|  \cJ|_{\overline{\Theta_n^{I}} H^{p'}(\CC_-, \HH) \rightarrow L^{s'}(\CC_+, \tilde \mu_{n,s}, \HH)}\| 
     \quad(n \in \NN),
$$
$$
 m_{p,s}
          =  \|  \cJ|_{\overline{\Theta^{I}} H^{p'}(\CC_-, \HH)
                             \rightarrow L^{s'}(\CC_+, \tilde \mu_{s}, \HH)}\| \quad(n \in \NN),
$$
where
$$
   \tilde \mu_{n,s} = \sum_{k=1}^n
     \frac{|2\re  z_k|^{s'}}{|b_{n,k}|^{s'}}  (G_k^{-1}(G_k^{-1})^*)^{s'/2}
                                                                      \delta_{z_k},
$$
$$
   \tilde \mu_{s} = \sum_{k=1}^\infty
     \frac{|2\re  z_k|^{s'}}{|b_{\infty,k}|^{s'}}  (G_k^{-1}(G_k^{-1})^*)^{s'/2}
                                                                      \delta_{z_k},
$$
\end{corollary}
\proof
This follows immediately from Theorem \ref{thm:carlesone}.
\qed

We have thus reduced the interpolation problem to the boundedness of an operator-weighted Carleson embedding.
In the finite-dimensional case,
we can in many instances give criteria for
the boundedness of this embedding:

\begin{theorem}   \label{thm:interpolation}
Let $1 \le s,p < \infty$, let $\HH = \CC^N$,
 let $(G_k)_{k \in \NN}$, $(I_k)_{k \in \NN}$, $(z_k)_{k \in \NN}$,
be defined as above, and let
$f \mapsto E(f) = (G_k f(z_k))_{k \in \NN}$ be the evaluation operator.
\begin{enumerate}
\item Let $1 \le p < \infty$, $1< s < \infty$.
Then $E( H^p(\CC_+, \CC^N)) \supseteq \ell^s(J_k)$, if and only if the
scalar measure
$$ \sum_{k=1}^\infty
     \frac{|2\re  z_k|^{s'}}{|b_{\infty,k}|^{s'}} \|(G_k^{-1})^* \overline{\Theta^I(z_k)}\|^{s'}\delta_{z_k}
$$
is an  $\frac{s'}{p'}$-Carleson measure. (This holds also in case
$p=1$ with the notion of $0$-Carleson measure from Equation (\ref{eq:zerocarl})).
\item
Let $s=1$, $1 \le p < \infty$. Then $E( H^p(\CC_+, \HH)) \supseteq \ell^1(J_k)$, if and
only if the operator sequence
$$
     \left(\frac{|\re  z_k|^{1/p}}{|b_{\infty,k}|}(G_k^{-1})^* \overline{\Theta^I(z_k)}\right)
$$
is bounded.
\end{enumerate}
\end{theorem}
\proof A weak$^*$ compactness argument shows that $E(H^p(\CC_+,
\HH)) \supseteq \ell^s(J_k)$ if and only if $(m_{n,p,s})$ is bounded. The
remainder of the first part follows directly from the comparison of the norm and the
trace of a positive matrix, Theorem \ref{thm:carlesone}, Theorem
\ref{thm:opcarleson} and the invariance of the $\alpha$-Carleson condition under
complex conjugation. For the second part, recall additionally from the
proof of Lemma \ref{lemma:carlfinite} that
\begin{multline*}
   m_{n,p,1}= \sup_{f \in H^{p'}(\CC_-, \HH) , \|f\|=1}
          \sup_{k \le n} \| 2 \frac{\re (z_k)}{|b_{n,k}|}  (G_k^{-1})^*
   \overline{\Theta^I_n(z_k)}  {f}(z_k)\| \\
= \sup_{k \le n} \| 2 \frac{\re (z_k)}{|b_{n,k}|}  (G_k^{-1})^*
   \overline{\Theta^I_n(z_k)}\| \|k_{z_k}\|_{L^p(i\RR)/H^p(\CC_-)} \approx \sup_{k \le n} \frac{|\re (z_k)|^{1/p}}{|b_{n,k}|}  \|(G_k^{-1})^*
   \overline{\Theta^I_n(z_k)}\|.
\end{multline*}
\qed

With Theorem
\ref{thm:opcarleson}, we further obtain
\begin{corollary}\label{cor:infdim}
If $\HH$ is a separable Hilbert space, $1 \le p < \infty$, $1< s< \infty$, and
$$
   \sum_{k=1}^\infty
     \delta_{z_k} \frac{|2\re  z_k|^{s'}\|{G_k^{-1}}^* \Theta_n^I(z_k)\|^{s'}}{|b_{\infty,k}|^{s'}} ,
$$
is a scalar $ \frac{s'}{p'}$-Carleson measure,
then $E( H^p(\CC_+, \CC^N)) \supseteq \ell^s(J_k)$.
\end{corollary}
\qed

If $p=2$ and $N = \dim \HH < \infty$, we can also deal with the weighted case.
\begin{theorem}   \label{thm:winterpolation}
Let $(G_k)_{k \in \NN}$, $(I_k)_{k \in \NN}$, $(z_k)_{k \in \NN}$,
be defined as above, and let
$f \mapsto E(f) = (G_k f(z_k))_{k \in \NN}$ be the evaluation operator.
Suppose that ${\Theta^{I^\perp}}^* W \Theta^{I^\perp}$ is a matrix $A_2$ weight.

Then $E( H^2_W(\CC_+, \CC^N)) \supseteq \ell^2(J_k)$, if and only if
the scalar measure
$$
    \sum_{k=1}^\infty
     \delta_{z_k}
         \frac{|2\re  z_k|^2}{|b_{\infty,k}|^2}
\|{G_k^{-1}}^* \Theta^I(z_k) ((\overline{{\Theta^{I^\perp}}^* W \Theta^{I^\perp}})(z_k))^{1/2}\|^2
$$
is Carleson.
\end{theorem}

\proof
By Theorem \ref{thm:carlesone}, we have to investigate the boundedness of the Carleson embedding
\begin{equation}   \label{eq:weightemb}
  \cJ_{\mu}: H_{\tilde W}^{2}(\CC_+, \CC^N) \rightarrow L^{2}(\CC_+, \CC^N,\mu),
\end{equation}
where
$
   \tilde W = \overline{ {\Theta^{I^\perp}}^* W^{-1} \Theta^{I^\perp} }
$
and
$
\mu =\sum_{k=1}^\infty
     \delta_{z_k}
         \frac{|2\re  z_k|^2}{|b_{\infty,k}|^2}
       {\Theta^I(z_k)}^* {G_k^{-1}} {G_k^{-1}}^* \Theta^I(z_k) .
$ The weight $\tilde W$ is a matrix $A_2$ weight, since
${\Theta^{I^\perp}}^* W \Theta^{I^\perp}$ is a matrix $A_2$ weight.
With the notation as in Theorem \ref{thm:weightcarleson}, we have
$$
  L^{2}(\CC_+, \CC^N,\mu)= L_{\tilde W}^{2}(\CC_+, \CC^N,\mu_{{\tilde W}^{-1}})
$$
where
$
d \mu_{{\tilde W}^{-1}}(z) = ({\tilde W}(z))^{-1/2} d \mu(z) ({\tilde W}(z))^{-1/2}.
$
Thus by Theorem \ref{thm:weightcarleson}, \ref{thm:opcarleson} and a comparison of trace and norm,
the embedding
\begin{equation}    \label{eq:l2emb}
 L_{\tilde W}^{2}(i \RR, \CC^N) \rightarrow L^{2}(\CC_+, \CC^N,\mu),
\end{equation}
is bounded if and only if the scalar measure
\begin{equation}    \label{eq:meas}
\sum_{k=1}^\infty
     \delta_{z_k}
         \frac{|2\re  z_k|^2}{|b_{\infty,k}|^2}
\|{G_k^{-1}}^* \Theta^I(z_k) (\tilde W (z_k))^{-1/2}\|^2
\end{equation}
is Carleson. By the splitting $L^2_{\tilde W}(i \RR, \CC^N) \simeq
H^2_{\tilde W}(\CC_+, \CC^N) \oplus H^2_{\tilde W}(\CC_-,
       \CC^N)$ as in (\ref{eq:tv}), the latter
embedding (\ref{eq:l2emb}) is bounded if and only if
(\ref{eq:weightemb}) is bounded.

Finally, using the $A_2$ property of $\tilde W$ again, we see that
the measure (\ref{eq:meas}) can be replaced by
\begin{multline*}
\sum_{k=1}^\infty
     \delta_{z_k}
         \frac{|2\re  z_k|^2}{|b_{\infty,k}|^2}
\|{G_k^{-1}}^* \Theta^I(z_k) (\tilde W^{-1} (z_k))^{1/2}\|^2 \\
=
\sum_{k=1}^\infty
     \delta_{z_k}
         \frac{|2\re  z_k|^2}{|b_{\infty,k}|^2}
\|{G_k^{-1}}^* \Theta^I(z_k) ((\overline{{\Theta^{I^\perp}}^* W \Theta^{I^\perp}})(z_k))^{1/2}\|^2.
\end{multline*}

\qed
 
In contrast to the scalar case, the matrix $A_2$ property for $W$ does not necessarily
imply the $A_2$ property for $ {\Theta^{I^\perp}}^* W \Theta^{I^\perp}$. In fact, it is not difficult
to show that if $W$ is matrix $A_2$, then $ {\Theta^{I^\perp}}^* W \Theta^{I^\perp}$ is
matrix $A_2$ if and only if the multiplication operator
$M_{\Theta^{I^\perp}}: L^2_W(i \RR) \rightarrow L^2_W(i \RR)$ is bounded. One easily sees that such
multiplication operators can have arbitrarily large norm even for weights of the form 
$$
   \left(   \begin{matrix} 1 & 0 \\ 0 & |\w|^\alpha \\   \end{matrix} \right)
$$ 
for fixed $\alpha$, $|\alpha| <1$.

Therefore, the most important case for applications is the case of scalar weights. Here we can also deal with
$1<p< \infty$, following \cite{mcphail}. We will state our condition in a slightly different way from
\cite{mcphail}. Recall that for $1<p<\infty$
a function $w$ on $\RR$ is called an \emph{$A_p$ weight},
if it is measurable, a.e. positive, locally integrable, and
$$
     \sup_{I \subset \RR, I \text{ interval }} 
               \left( \frac{1}{|I|} \int_I w(x) dx \right)  \left( \frac{1}{|I|} \int_I w^{-1/(p-1)}(x) dx \right) < \infty. 
$$
or, in M\"obius-invariant form, if
$$
  \sup_{z \in \CC_+} w(z) (w^{-1/(p-1)})(z) < \infty,
$$
where the symbols $w$, $w^{-1/(p-1)}$ are also used for the harmonic extensions of the respective weights
 to $\CC_+$.
It follows immediately from the definition that $w$ is an $A_p$-weight if and only if $w^{-1/(p-1)}$ is an $A_{p'}$-weight.
With this notation, we obtain

\begin{corollary}   \label{cor:wpinterpolation}
Let $(G_k)_{k \in \NN}$, $(I_k)_{k \in \NN}$, $(z_k)_{k \in \NN}$,
be defined as above, and let
$f \mapsto E(f) = (G_k f(z_k))_{k \in \NN}$ be the evaluation operator.
Suppose that $1<p<\infty$ and that $w$ is a scalar $A_p$ weight on $i \RR$.

Then $E( H^p_w(\CC_+, \CC^N)) \supseteq \ell^p(J_k)$, if and only if
the scalar measure
$$
    \sum_{k=1}^\infty
     \delta_{z_k}
         \frac{|2\re  z_k|^2}{|b_{\infty,k}|^{p'}}
\|{G_k^{-1}}^* \Theta^I(z_k) \|^{p'} w(z_k)
$$
is Carleson.
\end{corollary}
\proof
By Theorem \ref{thm:carlesone}, we have to investigate the boundedness of the Carleson embedding
\begin{equation}   \label{eq:pweightemb}
  \cJ_{\mu_p}: H_{w^{-1/(p-1)}}^{p'}(\CC_+, \CC^N) \rightarrow L^{p'}(\CC_+, \CC^N,\mu_p),
\end{equation}
where
$
\mu_p =\sum_{k=1}^\infty
     \delta_{z_k}
         \frac{|2\re  z_k|^{p'}}{|b_{\infty,k}|^{p'}}
       ({\Theta^I(z_k)}^* {G_k^{-1}} {G_k^{-1}}^* \Theta^I(z_k))^{p'/2} .
$
Writing $\tilde w$ for the harmonic extension of $w^{-1/(p-1)}$ and, similarly to the previous proof, 
$$
  L^p(\CC_+, \CC^N, d\mu_p)=   L^p_{\tilde w}(\CC_+, \CC^N, d\tilde \mu_p), 
$$
where
$ d \tilde \mu_p(z) = \tilde w^{-1}(z) d \mu_p(z)$, we can use the fact that $\tilde w$ is an $A_{p'}$-weight
together with the scalar weighted Carleson embedding
theorem to obtain that the embedding is bounded, if and only if $d\tilde \mu_p = \tilde w^{-1}(z) d \mu_p(z)$
is a Carleson measure. The M\"obius-invariant form of
the $A_p$ condition for $w$ ensures that $\tilde w^{-1}(z_k) \approx w(z_k)$ for all $k$,
 with equivalence constants only
depending on the $A_p$ constant of $w$, and we obtain that $d\tilde \mu_p$ is a Carleson measure if
and only if  
$\sum_{k=1}^\infty
     \delta_{z_k}
         \frac{|2\re  z_k|^2}{|b_{\infty,k}|^{p'}}
\|{G_k^{-1}}^* \Theta^I(z_k) \|^{p'} w(z_k)
$
is Carleson.
\qed

Returning to the case $p=2$, we obtain an interpolation result for certain Sobolev spaces which will be
useful for an application in control theory.

For $\beta >0$, recall the definition of the Sobolev space $\HH^2_\beta(\RR_+)$,
$$
 \HH^2_\beta(\RR_+) = \{ f \in L^2(\RR_+): |x|^\beta \hat f \in L^2(i\RR) \}. 
$$
This is a Hilbert space with the norm $\|f\|^2_{2, \beta} = \|\hat f\|_2^2 + \| |x|^\beta\hat f\|_2^2$.
Letting $\LL$ denote the \emph{Laplace transform},
$\LL f(z) = \hat f(z)=\int_0^\infty f(t) e^{-tz} dt$ for $z \in \CC_+$, we obtain

\begin{corollary} \label{cor:sobolevinter}
 Let $(G_k)_{k \in \NN}$, $(I_k)_{k \in \NN}$, $(z_k)_{k \in \NN}$,
be defined as above. Let $0<\beta <1/2$, and let $E$ be the evaluation operator
on $\HH^2_\beta(\RR_+)$ given by
$f \mapsto E(f) = (G_k \LL f(z_k))_{k \in \NN}$. 
Then $E( \HH^2_{\beta}(\RR_+, \CC^N)) \supseteq \ell^2(J_k)$, if and only if
the scalar measure
$$
    \sum_{k=1}^\infty
     \delta_{z_k}
         \frac{|2\re  z_k|^2}{|b_{\infty,k}|^{2}}
\|{G_k^{-1}}^* \Theta^I(z_k) \|^{2} |\w|^{2\beta}(z_k)
$$
is Carleson.
\end{corollary}
\proof
Clearly the Laplace transform defines an operator
$$
   \HH^2_{\beta}(\RR_+, \CC^N) \rightarrow H^2_{(1+ |\w|^{2 \beta})}(\CC_+, \CC^N)
$$
which is an isometric isomorphism up to an absolute constant.
It is well-known that the weight $(1 + |\w|^{2 \beta})$ is $A_2$ if and only if $|\beta | < 1/2$. Thus we
obtain the result from Corollary \ref{cor:wpinterpolation}.
\qed

\subsection{Some estimates for $m_{n,p,s}$}

We give some estimates for the interpolation constant $m_{n,p,s}$. The
proofs are similar to the case $p=2$ in \cite{jpp06}, Sections 2.4
and 2.5, so we just state the notation and results here.

\subsubsection{Finite union of Carleson sequences}
The first estimate concerns the case that $(z_k)$ is the union of
$K$ Carleson sequences. Here, for an estimate of the $m_{n,p,s}$ (up
to a constant), the $\Theta^I_n(z_k)$ can be replaced by a
Blaschke--Potapov product with at most $K$ factors,
$\Theta^{I,r}_n(z_k)$, where the factors correspond to the $z_j$ in
a suitably small hyperbolic $r$-neighbour\-hood of $z_k$.

\begin{corollary} \label{cor:carlun}
Let $(z_k)$ be the union of $K$ Carleson sequences and let $r>0$ be such that
each of the Carleson sequences is $r$-separated in the hyperbolic metric. For $k \in \NN$, define
$\Theta_{n,z_k,r}^I$ as the Blaschke--Potapov product associated to the
shift-invariant subspace
\[
\LL_{n,z_k,r} = \{f \in H^p(\CC_+, \CC^N), f(z_j) \in I_j \text{ for
all } z_j \text{ with } d(z_j,z_k) < r/2, j \le n \}.
\]
Then for $1 < p < \infty$,
$$
 m_{n,p,s} \approx \| \cJ_{\mu_{n,I,r}}|_{H^{p'}(\CC_+, \CC^N)\rightarrow L^{s'}(\CC_+, \HH, \mu_{n,I,r,s})} \|
           \quad(n \in \NN),
$$
where
$\mu_{n,I,r,s} = \sum_{k=1}^n
     \frac{|2\re  z_k|^{s'}}{|b_{n,k}|^{s'}}
     \|G_k^{-1}\Theta^{I}_{n,z_k,r}(z_k)\|^{s'}  \delta_{z_k}$
and $\cJ_{\mu_{n,I,r}}$ is the associated Carleson embedding.
\end{corollary}
\proof As in the case $p=2$ in \cite{jpp06}, Corollary 2.12. \qed

\subsubsection{Angles between subspaces} As in the case $p=2$,
interpolation conditions can be written in terms of angles between
certain subspaces of $H^2(\CC^N)$ rather than in terms of the inner
function $\Theta^I$. Recall that the angle between
two non-zero vectors $v_1,v_2$ in a Hilbert space $V$ is given by
\[   \angle(v_1,v_2):=\arccos \left(\frac{\langle v_1, v_2\rangle}{\|v_1\|\,\|v_2\|} \right),\] 
and the angle between two nontrivial subspaces $V_1$ and $V_2$ of
$V$ is defined as
\[    \angle(V_1,V_2):= \inf_{v_1\in V_1\backslash \{0\}, v_2\in V_2\backslash \{0\}}
      \angle(v_1,v_2) .\]
The angle between a vector and a subspace is defined analogously.\\

To recall some notation, for $n \in \NN$, $k=1, \dots, n$, we write
$$
\KK_{k, I} = ((b_k P_{I_k} + P_{{I_k}^\perp}))^\perp = k_{z_k} I_k,
$$
$$
\KK_{k, I,n}'=
\span\{ k_{z_j} I_j: j=1, \dots, n, j \neq k \} = ({\Theta'}_{k,n}^{I^\perp} H^2(\CC^N))^\perp
$$
and
$$
\KK_{k, I}'=
\overline{\span\{ k_{z_j} I_j: j \in \NN, j \neq k \}}= ({\Theta'}_{k}^{I^\perp} H^2(\CC^N))^\perp.
$$
where ${\Theta'}_{k,n}^{I^\perp}$ is the Blaschke--Potapov product as in Lemma
\ref{lemm:blaschke} corresponding
to $\{z_j, j=1, \dots, n, j \neq k \}$, and
${\Theta'}_{k}^{I^\perp}$ is the infinite Blaschke--Potapov product corresponding
to $\{z_j, j \in \NN, j \neq k \}$. We will state some interpolation results in terms of angles between 
such subspaces in $H^2(\CC_+, \CC^N)$.

\begin{corollary}  \label{thm:anglechar} Suppose that $N = \dim \HH < \infty$.
Suppose that there is a sequence of positive real numbers
$(\alpha_k)$ such that with the above notation, $G_k^* G_k =
\alpha_k^2 P_{I_k}$ for all $k \in \NN$. Then for $1 < p,s < \infty$,
$$
     m_{n,p,s} \approx \carl_{s'/p'}( \sum_{k=1}^n
         \frac{(\re  z_k)^{s'}}{|\alpha_k|^{s'} |\angle(\KK_{k, I},\KK'_{k, I,n})|^{s'}}
                             \delta_{z_k})^{1/s'}
$$
and
$$
    m_{p,s} = \sup_{n \in \NN}
    m_{n,p,s} \approx \carl_{s'/p'}( \sum_{k=1}^\infty
            \frac{(\re  z_k)^{s'}}{|\alpha_k|^{s'} |\angle(\KK_{k, I},\KK'_{k, I})|^{s'}}
                             \delta_{z_k})^{1/s'}
$$
with equivalence constant depending only on $N$, $p$, and $s$.
\end{corollary}

In the case of $(z_k)$ being the union of $K$ Carleson sequences,
Corollary \ref{cor:carlun} yields
\begin{corollary}  \label{cor:angle}
Let $1 < p < \infty$, $G_k^* G_k = \alpha_k^2 P_{I_k}$ for all $k$,
let $(z_k)$ be the union of $K$ Carleson sequences, and let $r>0$ be
such that each of the Carleson sequences is $r$-separated in the
hyperbolic metric. For $k \in \NN$, define $\KK'_{k,I,r,n}=\span\{
k_{z_j} I_j: j=1, \dots, n, j \neq k, d(z_j, z_k) < r/2 \}$. Then
$$
     m_{n,p,s} \approx
      \carl_{s'/p'}( \sum_{k=1}^n \frac{(\re  z_k)^{s'}}{|\alpha_k|^{s'} |\angle(\KK_{k, I},\KK'_{k, I,r,n})|^{s'}}
    \delta_{z_k})^{1/s'}
$$
and
$$
     m_{p,s}=\sup_{n \in \NN} m_{n,p,s} \approx
      \carl_{s'/p'}( \sum_{k=1}^\infty \frac{(\re  z_k)^{s'}}{|\alpha_k|^{s'} |\angle(\KK_{k, I},\KK'_{k, I,r})|^{s'}}
    \delta_{z_k})^{1/s'}
$$
with equivalence constant depending only on $N,r,K$. Here, we define
$\KK'_{k,I,r}=\span\{ k_{z_j} I_j: j \in \NN, j \neq k, d(z_j, z_k)
< r/2 \}$ and $|\angle(\KK_{k, I^\perp},\KK'_{k, I,r,n})|= \pi/2$,
if $\{ z_j: j=1, \dots, n, j \neq k, d(z_j, z_k) < r/2 \}= \empset$.
\end{corollary}

In the case that $G_k^*G_k$ is not the multiple of an orthogonal
projection, the $m_{n,p,s}$ can still be estimated in terms of angles
between subspaces in $H^2(\CC_+, \CC^N)$, albeit in a more technical way.
\begin{corollary}  \label{thm:genangle} Suppose that $N = \dim \HH < \infty$,
and suppose that
for each $k \in \NN$, the operator $G_k: \CC^N \rightarrow \CC^N$ is given as
$$
      G_k = \left(\begin{matrix}g_{k,1}^* \\ \vdots \\ g_{k,d_k}^*\\0\\ \vdots\\0     \end{matrix}   \right).
$$
Then for $1 < p < \infty$,
$$
     m_{n,p,s} \approx \carl_{s'/p'}( \sum_{k=1}^n \left(\sum_{i=1}^{d_k}
         \frac{(\re  z_k)^2 \| G_k^{-1} e_i\|^2}{ |\angle(k_{z_k} V^{G}_{k,i} ,\KK'_{k, I,n})|^2}\right)^{s'/2}
                             \delta_{z_k})^{1/s'}
$$
and
$$
    m_{p,s} = \sup_{n \in \NN}
    m_{n,p,s} \approx \carl_{s'/p'}( \sum_{k=1}^\infty \left(\sum_{i=1}^{d_k}
            \frac{(\re  z_k)^2 \| G_k^{-1} e_i\|^2}{ |\angle(k_{z_k} V^G_{k,i},\KK'_{k, I})|^2}\right)^{s'/2}
                             \delta_{z_k})^{1/s'}
$$
with equivalence constant depending only on $N$. Here,
$$
V^G_{k,i} =
\span\{\span\{g_{k,1}, \dots, g_{k, d_k}\}^\perp \cup
   \span\{g_{k,j}: 1 \le j \le d_k, j \neq i\} \} ^\perp
$$ for $1 \le i \le d_k$.
\end{corollary}

In an infinite-dimensional version, we only have an upper bound for
the $m_{p,s}$ from Theorem \ref{thm:opcarleson}.

\begin{corollary}  \label{cor:anglechar2}
With the notation as above, $\HH$ a separable Hilbert space, $1 < p
< \infty$,
$$
     m_{n,p,s} \lesssim \carl_{s'/p'}( \sum_{k=1}^n
         \frac{(\re  z_k)^{s'} \|G_k^{-1}\|^{s'}}{|\angle(\KK_{k, I},\KK'_{k, I,n})|^{s'}}
                             \delta_{z_k})^{1/s'}
$$
and
$$
    m_{p,s}   \lesssim  \carl_{s'/p'}( \sum_{k=1}^\infty
            \frac{(\re  z_k)^{s'} \|G_k^{-1}\|^{s'}}{ |\angle(\KK_{k, I},\KK'_{k, I})|^{s'}}
                             \delta_{z_k})^{1/s'}.
$$
\end{corollary}

Finally, we briefly want to comment on the boundedness of the evaluation operator.
\begin{theorem}  \label{thm:otherest}
Let $\HH=\mathbb C^N$ and for each $k\in \mathbb N$, let $G_k: \CC^N
\rightarrow \CC^N$ with the notation as above. 
\begin{enumerate}
\item  For $1\le p,s < \infty$, the following are equivalent
\begin{enumerate}
\item $E(H^p(\mathbb C_+, \mathbb C^N))\subset  \ell^s(J_k)$.
\item The measure $\sum_{k=1}^\infty \|G_k\|^s \delta_{z_k}$ is $s/p$-Carleson.
\end{enumerate}

\item
 For $1 < p,s < \infty$, the following are equivalent.
\begin{enumerate}
\item $E(H^p(\mathbb C_+, \mathbb C^N))= \ell^s(J_k)$.
\item  
\begin{enumerate}
\item the linear maps
  \begin{equation}\label{eq3}
      \frac{1}{(\re z_k)^{1/p}} G_k^*: J_k \rightarrow I_k
  \end{equation}
  are uniformly bounded above and below,

\item
 $(z_k)_{k \in \mathbb N}$ is the union of at most $N$ Carleson sequences, 
  and there exists a constant $r>0$ such that
      the systems $\{ I_k: z_k \in D_r(a)\}$ are uniformly Riesz in $\CC^N$ for all $a \in \CC_+$(in other words, the system 
$\{ k_{z_k} I_k \}_{k \in \NN}$ is unconditional in $H^p(\CC_+, \CC^N)$),

\item $\sum_{k=1}^\infty (\re z_k)^{s'/p'} \delta_{z_k} $ is an $s'/p'$-Carleson measure \\
and $\sum_{k=1}^\infty (\re z_k)^{s/p} \delta_{z_k} $ is an $s/p$-Carleson measure.

\end{enumerate}
(The last condition is redundant in the case $p=s$).
\end{enumerate}

\end{enumerate}
\end{theorem}
\proof
 1. This follows easily from
 $$
     \|E(f)\|^s = \sum_{k=1}^\infty \|G_k f(z_k)\|^s =
     \int_{\CC_{+}} \langle d \mu^{2/s} f(z), f(z)^{s/2} \rangle,
 $$
 where $\mu = \sum_{k=1}^\infty (G_k^* G_k)^{s/2} \delta_{z_k}$, the
 matrix Carleson--Duren embedding Theorem \ref{thm:opcarleson}, and a
 comparison of trace and norm.

 2. We first show (a) $\Rightarrow$ (b). If $E:H^p(\mathbb
C_+, \mathbb C^N) \rightarrow \ell^s(J_k)$ is bounded and
surjective, then
$$E^*: \ell^{s'}(J_k) \rightarrow H^{p'}(\mathbb C_+,
\mathbb C^N), \quad (x_k) \mapsto \sum_{k \in \NN} k_{z_k} G_k^*x_k
$$
is bounded and bounded below. Applying $E^*$ to $(0, \dots, 0, x_k, 0, \dots)$, we see that
$\| k_{z_k} G_k^* x_k \|_{p'} \approx \|x_k\|$ for all $ k \in \NN$, $x_k \in I_k$ and that 
the linear maps
  $
      \frac{1}{(\re z_k)^{1/p}} G_k^*: J_k \rightarrow I_k
  $
  are uniformly bounded above and below. In other words, the map
  $$
    \ell^s(J_k) \rightarrow \ell^s(I_k),  \quad (x_k) \mapsto \left(\frac{1}{(\re z_k)^{1/p}}
  G_k^*x_k\right)
  $$
  is an isomorphism of Banach spaces. That means, the map
$$
   \tilde E: H^p(\mathbb C_+, \mathbb C^N))= \ell^s (I_k), \quad f \mapsto (\re z_k)^{1/p} P_{I_k} f(z_k)
$$
is bounded and surjective, and
\begin{equation}
   \tilde E^*: \ell^{s'}(I_k) \rightarrow H^{p'}(\CC_+, \CC^N), 
                \quad (x_k) \mapsto \sum_{k=1}^\infty (\re z_k)^{1/p} k_{z_k} x_k 
\end{equation}
is bounded and bounded below.

Boundedness of $\tilde E$ implies with Part (1) that the measure
$ \mu_1 = \sum_{k=1}^\infty (\re z_k)^{s/p} \delta_{z_k}$ is $s/p$-Carleson. Surjectivity of $\tilde E$
implies by Corollary \ref{thm:anglechar} that
$$
\sum_{k=1}^\infty
            \frac{(\re  z_k)^{s'}}{(\re z_k)^{s'/p} |\angle(\KK_{k, I},\KK'_{k, I})|^{s'}}
                             \delta_{z_k}
   =\sum_{k=1}^\infty
            \frac{(\re  z_k)^{s'/p'}}{ |\angle(\KK_{k, I},\KK'_{k, I})|^{s'}}
                             \delta_{z_k}
$$
is a $s'/p'$ Carleson measure, so certainly
$$
   \mu_2=\sum_{k=1}^\infty (\re  z_k)^{s'/p'} \delta_{z_k}
$$
is a $s'/p'$ Carleson measure. 
(Here, $\alpha_k = (\re z_k)^{1/p}$ in the notation of Corollary \ref{thm:anglechar}.)

The necessary condition (3) in Theorem \ref{thm:ccarleson} for both $\mu_1$ and $\mu_2$
and a simple convexity argument then imply that 
$\sum_{k=1}^\infty \re  z_k  \delta_{z_k}$ is a Carleson measure, and $(z_k)$ is consequently
a finite union of Carleson sequences (see e.g. \cite{Nik86}, Lecture VII).

Boundedness below of $\tilde E^*$ then implies that for suitable $r >0$,
the systems $\{ I_k : z_k \in D_r(a) \}$ are uniformly Riesz in $\CC^N$ for all $a \in \CC_+$.

By \cite{Tre86}, this means that the system $\{ k_{z_k} I_k \}_{k \in \NN}$ is unconditional in $H^2(\CC_+, \CC^N)$.

For the direction (b) $\Rightarrow$ (a), note that the unconditionality of
$\{ k_{z_k} I_k \}_{k \in \NN}$ in $H^2(\CC_+, \CC^N)$ implies in
particular
$$ 
     \inf_{k \in \NN} \angle(\KK'_{k,I}, \KK_{k,I}) > 0.
$$
So, since $\mu_1 = \sum_{k=1}^\infty (\re  z_k)^{s'/p'} \delta_{z_k}
$
is a $s'/p'$ Carleson measure, the measure
$$
\sum_{k=1}^\infty
            \frac{(\re  z_k)^{s'/p'}}{ |\angle(\KK_{k, I},\KK'_{k, I})|^{s'}}
                             \delta_{z_k}
$$
is also $s'/p'$-Carleson, and $\tilde E$ is surjective by Corollary \ref{thm:anglechar}.

Since $\mu_1 =\sum_{k=1}^\infty (\re z_k)^{s/p} \delta_{z_k}$ is $s/p$-Carleson, it follows from Part(1)
that $\tilde E$ is also bounded. The uniform boundedness and boundedness below of the
maps $\frac{1}{(\re z_k)^{1/p}} G_k^*: J_k \rightarrow I_k$ now imply boundedness and surjectivity of $E$.
\qed

\section{Controllability}
\label{sec:3}

In this section we apply the results on interpolation by vector-valued analytic functions to controllability problems of infinite-dimensional linear systems.
We study a system of the form
\begin{eqnarray}\label{system}
 \dot{x}(t) &=& Ax(t)+Bu(t), \qquad t\ge 0,\\
x(0) &=& x_0.\nonumber
\end{eqnarray}
Here we assume that $A$ is the generator of an exponentially stable
$C_0$-semigroup $(T(t))_{t\ge 0}$ on a Banach space
$X$ such that for some $s$ with $1 \le s < \infty$ 
the eigenvectors $(\phi_n)_{n\in\mathbb N}$ of $A$ form a basis of $X$, equivalent to the standard basis of $\ell^s$,
and the corresponding eigenvalues $(\lambda_n)_{n\in \mathbb N}$ are pairwise distinct. The eigenvalues
$(\lambda_n)_{n\in\mathbb N}$ then lie in the open left half plane uniformly bounded away from the imaginary axis.
For our input space $\cU$ we shall fix $\cU=L^p(0,\infty;\mathbb C^N)$ for some $p$ with $1<p<\infty$ or
a Sobolev space $\cU=\HH^2_\beta(\RR_+)$ with $-\frac12<\beta<\frac12$; then we take $u \in \cU$.
We assume that the control operator $B$ is given by
\[ Bv = \sum_{n=1}^\infty \langle v,b_n\rangle \phi_n, \qquad v\in \CC^N, \]
where  $(b_n)_n\subseteq \CC^N$, and, to avoid trivial cases, that $b_n \ne 0$ for all $n$.
Thus $B$ is a linear bounded operator from $\mathbb C^N$ to
\[ X_{(b_n)} =\left\{\sum_{n\in\mathbb N} x_n \phi_n:
  \left\{\frac{x_n|b_n|^{-1}}{n^{2/s}}\right\}_{n \in \NN} \in \ell^s\right\},\]
equipped with the norm
\[ 
\|x\|_{(b_n)} := \left(\sum_{n\in\mathbb N} \frac{|x_n|^s|b_n|^{- s}}{n^2} \right)^{1/s}.
\]

One important feature of the interpolation space $X_{(b_n)}$ is that
the semigroup $(T(t))_{t\ge 0}$ can be extended to a $C_0$-semigroup on
$X_{(b_n)}$, which we denote again by $(T(t))_{t\ge 0}$, using the property that $T(t)\phi_n=e^{\lambda_n t}\phi_n$ for $n \in \NN$,
and the generator of this
extended semigroup, denoted again by $A$, is an extension of $A$.
By a solution of the system (\ref{system}) we mean the so-called mild solution given by
\[ x(t)= T(t)x_0 +\int_0^t T(t-s)Bu(s)\,ds,\]
which is a continuous function with values in the interpolation space $X_{(b_n)}$.
We introduce
the operator
${\cal B}_{\infty}\in{\cal L}(\cU, X_{(b_n)})$ by
\[ {\cal B}_{\infty}u := \int_0^\infty T(s)B u(s)\,ds.\]
In the literature on infinite-dimensional system it is often assumed that
the operator $B$ is admissible for the semigroup $(T(t))_{t\ge 0}$,
and thus for some of our results we will include admissibility in
the assumptions.

\begin{definition}
$B$ is called {\em  admissible for $(T(t))_{t\ge 0}$}, if  ${\cal B}_{\infty} u \in H$ for every
$u\in \cU$.
\end{definition}

Admissibility implies that the mild solution of (\ref{system}) corresponding
to an initial condition $x(0)=x_0\in H$ and to $u\in \cU$ is a
continuous $H$-valued function of $t$. The case $\cU=L^p(0,\infty;\mathbb C^N)$ has been introduced and
studied in \cite{wei89} and \cite{un07}. The case  $\cU=\HH^2_\beta(\RR_+)$ seems to be new and not yet studied in the literature. For further information on admissibility
we refer the reader to the survey \cite{jp}. 

Using the special representation of $A$ and $B$ we see that
\begin{equation}\label{eqninterpolation}
\int_0^\infty T(t)Bu(t)\,dt
= \sum_{n=1}^\infty \int_0^\infty e^{\lambda_n t}\langle u(t), b_n\rangle dt\,\phi_n
= \sum_{n=1}^\infty \langle \hat{u}(-\lambda_n), b_n\rangle \phi_n,
\end{equation}
for every $u\in\cU$.

It follows that $B$ is admissible  for $(T(t))_{t\ge 0}$
if and only if
$ \wB \cU \subseteq \ell^s(\mathbb N)$,
where $\wB: \cU\rightarrow \{x:\mathbb N\rightarrow
\mathbb C\}$ is defined by
\begin{equation}\label{defu}
 \wB g := (\langle \hat g(-\lambda_n), b_n\rangle)_n.
 \end{equation}

We shall write $\LL\cU$ for the space of Laplace transforms of $\cU$, noting that for
$1 < p \le 2$ we have a bounded operator $\LL: L^p(0,\infty; \CC^N) \mapsto H^{p'}(\CC_+,\CC^N)$, where $p'$ is the conjugate index to $p$,
and for $2 \le p < \infty$ we have a bounded operator
$\LL^{-1}: H^{p'}(\CC_+,\CC^N) \to L^p(0,\infty; \CC^N)$. 
(This is basically the Hausdorff--Young inequality \cite[VI.3]{katznelson}.)
The other choice of $\cU$ mentioned above is
 the Sobolev space $\HH^2_\beta(\RR_+)$, and here
we have already noted
that $\LL$ defines an operator
$$
   \HH^2_{\beta}(\RR_+, \CC^N) \rightarrow H^2_{(1+ |\w|^{2 \beta})}(\CC_+, \CC^N)
$$
which is an isometric isomorphism.

We now apply the results on interpolation by vector-valued analytic functions to admissibility. To obtain necessary and sufficient conditions, we 
work first with the spaces $\cU=\LL^{-1}H^{p'}(\CC_+,\CC^N)$ with the norm induced from $H^{p'}$, and then with the spaces $\cU=L^p(0,\infty;\mathbb C^N)$.

\begin{theorem}\label{theounteregge}
Suppose that $\cU=\LL^{-1}H^{p'}(\CC_+,\CC^N)$ with $1<p<\infty$. Then $B$ is admissible  for $(T(t))_{t\ge 0}$ if and only if 
the measure $\sum_{k=1}^\infty |b_k|^s \delta_{z_k}$ is $s/p'$-Carleson.

If $\cU=L^p(0,\infty;\mathbb C^N)$ with $1<p\le 2$ then $B$ is admissible if  the measure $\sum_{k=1}^\infty |b_k|^s \delta_{-\lambda_k}$ is $s/p'$-Carleson.

If $\cU=L^p(0,\infty;\mathbb C^N)$ with $2\le p<\infty $ then the admissibility of $B$ implies that the  measure $\sum_{k=1}^\infty |b_k|^s \delta_{-\lambda_k}$ is 
$s/p'$-Carleson.
\end{theorem}

\beginpf
Choosing ${\cal H}:=\mathbb C^N$ and defining $G_k\in \mathbb C^{N\times N}$ by $G_k^*:=(b_k \,\, 0 \cdots 0)$, $k\in\mathbb N$, 
the theorem follows immediately from Theorem \ref{thm:otherest}.

\kasten

Theorem \ref{theounteregge} can also be found in \cite{un07}. We shall discuss the
following controllability concepts.

\begin{definition}
Let $\tau>0$. We say that the system (\ref{system}) is
\begin{enumerate}
\item {\em null-controllable in time $\tau$}, if $R( T(\tau))
  \subset R({\cal B}_{\infty})$;
\item {\em approximately controllable}, if $R({\cal
    B}_{\infty}) \cap X$ is dense in $X$;
\item {\em  exactly controllable}, if $X\subset R({\cal
    B}_{\infty})$.
\end{enumerate}
\end{definition}

Here $R(\cdot)$ denotes the range of an operator.
It is easy to see that every exactly controllable system is
approximately controllable and null-controllable in any time $\tau>0$.

\subsection{Conditions for exact controllability}\label{subsec1}

As in \cite{jp05,jpp06}
we may reduce the question of exact controllability to an interpolation problem. This can then be solved using the results
of Section \ref{sec:2}.
Using (\ref{eqninterpolation}), it follows that
the system (\ref{system}) is  exactly controllable if and only if
$ \ell^s(\mathbb N) \subseteq \wB \cU$.
where $\wB: \cU\rightarrow \{x:\mathbb N\rightarrow
\mathbb C\}$ is defined by (\ref{defu}).

To obtain necessary and sufficient conditions, we 
work first with the input spaces $\cU=\LL^{-1}H^{p'}(\CC_+,\CC^N)$ (with the norm induced from $H^{p'}$). There are two cases to consider.

\begin{theorem}\label{theo2a}
Suppose that $\cU=\LL^{-1}H^{p'}(\CC_+,\CC^N)$ with $1< p \le s' < \infty$.
Then the following statements are equivalent:
\begin{enumerate}
\item System (\ref{system}) is exactly controllable.
\item There exists a constant $m>0$ such that for all $h>0$ and all $
  \omega\in\mathbb R$:
\begin{equation}\label{condition1a}
 \sum_{-\lambda_n\in R(\omega,h)}
\frac{|\re  \lambda_n|^{s'}}{\|b_n\|^{s'} |\angle( e^{\lambda_n t} b_n, {\rm span}_{j\not=n, j\in\mathbb N} \{ e^{\lambda_j t} b_j\} )|^{s'}} \le m h^{s'/p},
 \end{equation}
where $R(\omega,h):= \{s\in\mathbb C_+: \mbox{\rm Re}\,s <h,\
\omega-h<\mbox{\rm Im}\,s< \omega+h\}$.
\end{enumerate}
\end{theorem}

\begin{theorem}\label{theo2b}
Suppose that $\cU=\LL^{-1}H^{p'}(\CC_+,\CC^N)$ with $1 < s' < p < \infty$.
Then the following statements are equivalent:
\begin{enumerate}
\item System (\ref{system}) is exactly controllable.
\item The function
\begin{equation}\label{condition1b}
\omega \mapsto  \sum_{n \in \NN} 
\frac{|\re  \lambda_n|^{s'}}{\|b_n\|^{s'} |\angle( e^{\lambda_n t} b_n, {\rm span}_{j\not=n, j\in\mathbb N} \{ e^{\lambda_j t} b_j\} )|^{s'}} p_{-\lambda_n}(i\omega)
 \end{equation}
lies in $L^{p/(p-s')}(\RR)$. Here $p_{-\lambda_n}$ denotes the Poisson kernel (\ref{poisson}) for $-\lambda_n$.
\end{enumerate}
\end{theorem}

\begin{remark}{\rm In the scalar case $N=1$, expressions (\ref{condition1a}) and (\ref{condition1b}) can be simplified, since
\begin{equation}\label{eq:angleequiv}
 |\angle( e^{\lambda_n t} b_n, {\rm span}_{j\not=n, j\in\mathbb N} \{ e^{\lambda_j t} b_j\} )| \asymp 
\prod_{j \ne n} \left| \frac{\lambda_n-\lambda_j}{\overline{\lambda_n}+\lambda_j}\right|.
\end{equation}
The resulting expressions provide a generalization of \cite[Thm.~3.1]{jp05}.}
\end{remark}

{\bf Proof of Theorems \ref{theo2a} and \ref{theo2b}}
We choose ${\cal H}:=\mathbb C^N$ and we define $G_k\in \mathbb C^{N\times N}$ by $G_k^*:=(b_k \,\, 0 \cdots 0)$, $k\in\mathbb N$.
Note that system (\ref{system}) is exactly controllable if and only if $\ell^s(\mathbb N)\subset \wB \cU$. A weak$^*$ compactness argument shows 
that the latter holds if and only if
\[ \sup_{n\in\mathbb N}
\sup_{\substack{x\in\ell^s(\mathbb C^N)\\ \|x\|_s\le 1}} \inf \left\{ \|f\|_\cU : f\in \cU, G_k\hat f(-\lambda_k) = (x_k \,\,0 \cdots 0)^T, k=1,\cdots,n\right\}\]
is finite.
Thus we have reduced the question of exact controllability to an interpolation problem treated in Section 2.
Using the notation of Section 2 we have
\[ \angle({\cal K}_{k,I}, {\cal K}_{k,I}')=\angle(e^{\lambda_k t}b_k, \span_{j\not=k, j\in\mathbb N}\{ e^{\lambda_j t}b_j\} ).\]
Theorems \ref{theo2a} and \ref{theo2b} now follow from Corollary \ref{thm:anglechar}. 
\kasten

This gives an immediate corollary for $\cU=L^p(0,\infty; \CC^N)$. In the case $p=2$ it provides necessary and sufficient conditions for
controllability, but even for other values of $p$ it provides an implication in one direction or the other.

\begin{corollary}
Suppose that $\cU=L^p(0,\infty; \CC^N)$. Then:\\
(i) If $1<p \le s'<\infty$ and $p\ge 2$, then (\ref{condition1a}) is a sufficient condition for the exact controllability of (\ref{system}).\\
(ii) If $1<p \le s' < \infty$ and $p \le 2$, then (\ref{condition1a}) is a necessary condition for the exact controllability of (\ref{system}).\\
(iii) If $1 < s' < p < \infty$ and $p \ge 2$, then (\ref{condition1b}) is a sufficient condition for the exact controllability of (\ref{system}).\\
(iv) If $1 < s' < p < \infty$ and $p \le 2$, then  (\ref{condition1b}) is a necessary condition for the exact controllability of (\ref{system}).
\end{corollary}
\beginpf 
This follows from Theorems \ref{theo2a} and \ref{theo2b}, together with the Hausdorff--Young theorem mentioned above.
\endpf

We now consider controllability in the situation $\cU=\HH^2_{\beta}(\RR_+, \CC^N)$ (again sufficient conditions and necessary conditions can
be derived for other values of $p$ with $1<p<\infty$, using Corollary  \ref{cor:wpinterpolation}, but we shall omit them). 
For $0<\beta<1/2$ we shall write $\phi_{2\beta}(z)$ for the harmonic extension to $\CC_+$ of the function $i\omega \mapsto |\omega|^{2\beta}$
defined on the imaginary axis. This is given by
\[
\phi_{2\beta}(x+iy)=\frac{1}{\pi} \int_{-\infty}^\infty \frac{x |\omega|^{2\beta}}{(y-\omega)^2+x^2} \, d\omega.
\]
For real arguments (i.e., $y=0$), we have the particularly simple expression
\[
\phi_{2\beta}(x)= 
\frac{x^{2\beta}}{\pi}
\int_{-\infty}^\infty \frac{|t|^{2\beta}}{t^2+1} \, dt,
\]
which makes the analysis of the case of real eigenvalues by means of the following result relatively straightforward.

\begin{theorem}\label{sobolev}
Suppose that $\cU=\HH^2_{\beta}(\RR_+, \CC^N)$ with $0<\beta<1/2$.
Then the following statements are equivalent:
\begin{enumerate}
\item System (\ref{system}) is exactly controllable.
\item There exists a constant $m>0$ such that for all $h>0$ and all $
  \omega\in\mathbb R$:
\begin{equation}\label{conditionsob}
 \sum_{-\lambda_n\in R(\omega,h)}
\frac
{|\re  \lambda_n|^{2} \phi_{2\beta}(-\lambda_n)}
{\|b_n\|^{2} |\angle( e^{\lambda_n t} b_n, {\rm span}_{j\not=n, j\in\mathbb N} \{ e^{\lambda_j t} b_j\} )|^{2}}
 \le m h,
 \end{equation}
where $R(\omega,h):= \{s\in\mathbb C_+: \mbox{\rm Re}\,s <h,\
\omega-h<\mbox{\rm Im}\,s< \omega+h\}$.
\end{enumerate}
\end{theorem}

\beginpf
This is proved in the same way, using Corollary \ref{cor:sobolevinter} and the identity 
\[
\frac{\|G_n^{-1}\Theta^I(-\lambda_n)\|^2}{|b_{\infty n}|^2}
= \frac{1}{\|b_n\|^2} \frac{\|\Theta^I(-\lambda_n)\|^2}{|b_{\infty n}|^2}
\asymp \frac{1}{\|b_n\|^2}\frac{1}{ |\angle( e^{\lambda_n t} b_n, {\rm span}_{j\not=n, j\in\mathbb N} \{ e^{\lambda_j t} b_j\} )|^{2}},
\]
given in \cite[Lem~2.15]{jpp06}.
\endpf

Concerning admissibility and controllability we have the following equivalent condition.

\begin{theorem}
Suppose that $\cU=\LL^{-1}H^{p'}(\CC_+,\CC^N)$.
Then $B$ is an admissible control operator  and the system (\ref{system}) is exactly controllable if and only if
\begin{enumerate}
\item The sequence 
\[ \frac{\|b_k\|}{|{\rm Re}\, \lambda_k|^{1/p}}, \quad k\in\mathbb N,\]
is uniformly bounded above and below.
\item $\{k_{-\lambda_k} b_k\}$ is unconditional in $H^p(\mathbb C_+,\mathbb C^N)$.
\item $\sum_{k=1}^\infty (\re -\lambda_k)^{s'/p'} \delta_{-\lambda_k} $ is an $s'/p'$-Carleson measure \\
and $\sum_{k=1}^\infty (\re -\lambda_k)^{s/p} \delta_{-\lambda_k} $ is an $s/p$-Carleson measure.
\end{enumerate} 
\end{theorem}
\beginpf
Choose ${\cal H}:=\mathbb C^N$ and defining $G_k\in \mathbb C^{N\times N}$ by $G_k^*:=(b_k \,\, 0 \cdots 0)$, $k\in\mathbb N$, 
the theorem follows immediately from Theorem \ref{thm:otherest}.
\endpf

\subsection{Conditions for null controllability}
\label{sec:3.2}

As for exact controllability, the question of null controllability is easily reduced to an interpolation problem.
Using (\ref{eqninterpolation}) it is easy to see that the system (\ref{system}) is null-controllable in time $\tau$ if and only if 
$\{(e^{\lambda \tau}x_n)_n: (x_n)_n\in\ell^s(\mathbb N)\}\subset \wB \cU$, where $\wB$ is defined by (\ref{defu}).
Replacing $b_k$ by $e^{-\overline{\lambda_k} \tau}b_k$ in the previous subsection, we obtain the following two theorems (\ref{thm:3.8new} and \ref{thm:3.9new}).

\begin{theorem}\label{thm:3.8new}
Suppose that $\cU=\LL^{-1}H^{p'}(\CC_+,\CC^N)$ with $1< p \le s' < \infty$.
Then the following statements are equivalent:
\begin{enumerate}
\item System (\ref{system}) is null-controllable in time $\tau$.
\item There exists a constant $m>0$ such that for all $h>0$ and all $
  \omega\in\mathbb R$:
\begin{equation}\label{condition1an}
 \sum_{-\lambda_n\in R(\omega,h)}
\frac{|\re  \lambda_n|^{s'}
e^{s'{\rm Re}\,\lambda_n \tau}}
{\|b_n\|^{s'} |\angle( e^{\lambda_n t} b_n, {\rm span}_{j\not=n, j\in\mathbb N} \{ e^{\lambda_j t} b_j\} )|^{s'}} \le m h^{s'/p},
 \end{equation}
where $R(\omega,h):= \{s\in\mathbb C_+: \mbox{\rm Re}\,s <h,\
\omega-h<\mbox{\rm Im}\,s< \omega+h\}$.
\end{enumerate}
\end{theorem}

\begin{theorem}\label{thm:3.9new}
Suppose that $\cU=\LL^{-1}H^{p'}(\CC_+,\CC^N)$ with $1 < s' < p < \infty$.
Then the following statements are equivalent:
\begin{enumerate}
\item System (\ref{system}) is null-controllable in time $\tau$.
\item The function
\begin{equation}\label{condition1bn}
\omega \mapsto  \sum_{n \in \NN} 
\frac{|\re  \lambda_n|^{s'}e^{s'{\rm Re}\,\lambda_n \tau}}
{\|b_n\|^{s'} |\angle( e^{\lambda_n t} b_n, {\rm span}_{j\not=n, j\in\mathbb N} \{ e^{\lambda_j t} b_j\} )|^{s'}} p_{-\lambda_n}(i\omega)
 \end{equation}
lies in $L^{p/(p-s')}(\RR)$.
\end{enumerate}
\end{theorem}

\begin{remark}\label{rem:312}
{\rm Once again, expressions (\ref{condition1an}) and (\ref{condition1bn}) simplify
in the scalar case $N=1$, 
using (\ref{eq:angleequiv}); the resulting formulae provide  a generalization of \cite[Thm.~2.1]{jp05}.}
\end{remark}

\begin{corollary}\label{cor:313}
Suppose that $\cU=L^p(0,\infty; \CC^N)$. Then:\\
(i) If $1<p \le s'<\infty$ and $p\ge 2$, then (\ref{condition1an}) is a sufficient condition for the null controllability of (\ref{system}) in time $\tau$.\\
(ii) If $1<p \le s' < \infty$ and $p \le 2$, then (\ref{condition1an}) is a necessary condition for the null controllability of (\ref{system}) in time $\tau$.\\
(iii) If $1 < s' < p < \infty$ and $p \ge 2$, then (\ref{condition1bn}) is a sufficient condition for the null controllability of (\ref{system}) in time $\tau$.\\
(iv) If $1 < s' < p < \infty$ and $p \le 2$, then  (\ref{condition1bn}) is a necessary condition for the null controllability of (\ref{system}) in time $\tau$.
\end{corollary}
\beginpf 
Again this follows from Theorems \ref{thm:3.8new} and \ref{thm:3.9new}, using the Hausdorff--Young theorem.
\endpf

Similarly to Theorem \ref{sobolev}, equivalent conditions concerning null controllability for the case $\cU=\HH^2_{\beta}(\RR_+, \CC^N)$ with $0<\beta<1/2$ can be obtained.

\subsection{Conditions for approximate controllability}

Next we characterize approximately controllable systems in terms of their
eigenvalues and the operator $B$. By $e_n$ we denote the $n$th unit vector of $\mathbb C^N$. For the purposes of this
subsection, we introduce the interpolation space $X_{s,\alpha}$ defined for $\alpha \in \RR$ and $1 < s < \infty$ by
\[
X_{s,\alpha} = \left\{ \sum_{n \in \NN} x_n \phi_n: \{x_n |\lambda_n|^\alpha \} \in \ell^s \right\},
\]
with norm 
\[
\|x\|_{s,\alpha} := \left( \sum_{n\in \NN} |x_n|^s|\lambda_n|^{\alpha s} \right)^{1/s}.
\]
The dual space to $X_{s,\alpha}$, with the natural pairing, can be identified with $X_{s',-\alpha}$, and clearly $X=X_{s,0}$.
\begin{theorem}
Suppose that $\cU=L^p(0,\infty; \CC^N)$ with $1<p<\infty$, 
$\overline{\{\lambda_n: n\in \mathbb N\}}$ is {\em totally disconnected}, that is, no 
two points $\lambda, \mu\in \overline{\{\lambda_n: n\in \mathbb N\}}$ can be joined by a segment lying 
entirely in $\overline{\{\lambda_n: n\in \mathbb N\}}$. 
Then for $B \in {\cal L}(\mathbb
C^N,X_{s,\alpha})$ the following properties are equivalent:
\begin{enumerate}
\item The system (\ref{system}) is approximately controllable.
\item $\rank (\langle Be_1,\phi_n\rangle, \cdots, \langle Be_N,\phi_n\rangle)=1$ for all $n\in\mathbb N$.
\end{enumerate}
\end{theorem}

\beginpf It is easy to see that statement 1 implies statement 2.

To show that statement 2 implies statement 1, we adapt the proof of \cite[Thm.~4.2.3]{cuzw95}, beginning with the 
special case that $B \in {\cal L}(\mathbb C^N,X)$. We need to show that the reachability subspace 
$\cR=R({\cal
    B}_{\infty})$ is dense in $X$. As in \cite[Thm.~4.1.19]{cuzw95} we obtain that $\overline{\cR}$ is the smallest closed, 
$T(t)$-invariant subspace in $X$ containing $R(B)$ and hence equal to
the closed linear span of $\{\phi_n: n \in \JJ\}$ for some $\JJ \subseteq \NN$, see \cite[Thm.~2.5.8]{cuzw95}.
The remainder of the proof follows exactly as in \cite[Thm.~4.2.3]{cuzw95}.

To deduce the result in the general case $B \in {\cal L}(\mathbb
C^N,X_{s,\alpha})$, we fix an integer $m>-\alpha$.
Now we know that the system $(A, \beta)$, where
\[\beta:=\sum_{j=1}^\infty \frac{\langle \cdot,B^* \phi_j\rangle}{(1-\lambda_j)^m}  \phi_j\in {\cal L}(\mathbb C^N, X),\]
 is approximately controllable, by the
arguments above. Using the fact that
\[
 {\cal B}_\infty f =\sum_{j=1}^\infty    \int_0^\infty
e^{\lambda_n  t}
  \langle f(t), B^*\phi_j \rangle \, dt\,\phi_j = \sum_{j=1}^\infty
 \langle \hat{f}(-\lambda_j
  ),B^*\phi_j \rangle \phi_j,
\]
where $\hat{f}$ denotes the Laplace transform of $f$, we get that
the set
\[
S_\beta:=\left\{\sum_{j=1}^\infty \frac{\langle \hat{f}(-\lambda_j
  ),B^* \phi_j\rangle}{(1-\lambda_j)^m} \phi_j: \, f \in L^p(\RR_+,\mathbb C^N)\right\}
\]
is dense in $H$. Similarly, let
\[
S_B:=\left\{\sum_{j=1}^\infty \langle \hat{f}(-\lambda_j
  ),B^* \phi_j\rangle \phi_j: \, f \in L^p(\RR_+,\mathbb C^N)\right\}.
\]
Now
 if $f \in L^p(0,\infty; \mathbb C^N)$, then so is the function $g$ obtaining by taking the convolution
between $f$ and the function $t \mapsto t^{m-1}e^{-t}/(m-1)!$,
and then
\[
\sum_{j=1}^\infty \frac{\langle \hat{f}(-\lambda_j
  ),B^* \phi_j\rangle}{(1-\lambda_j)^m} \phi_j = \sum_{j=1}^\infty \langle \hat{g}(-\lambda_j
  ),B^* \phi_j\rangle \phi_j.
\]
Hence $S_\beta \subseteq S_B$, which implies that $S_B$ is dense in $X$, as required.
\kasten

It would be of interest to decide whether the above result still holds for arbitrary $X_{(b_n)}$, but it seems
that the methods of the proof do not extend directly to the general situation. 

\subsection{Application to the Heat Equation}

As in \cite{jp05,jpp06}, we shall very briefly consider the one-dimensional heat equation
on $[0,1]$, given by 
\begin{eqnarray*}
\frac{\partial z}{\partial t}(\xi,t)&=&\frac{\partial^2 z}{\partial \xi^2}(\xi,t), \qquad (\xi \in (0,1), \quad t \ge 0),\\
z(0,t)&=&0, \qquad z(1,t)=u(t), \qquad (t \ge 0),\\
z(\xi,0)&=& z_0(\xi), \qquad (\xi \in (0,1)).
\end{eqnarray*}
This may be written in the form (\ref{system}) with $X=L^2(0,1)$, and 
$A\phi_n=\lambda_n \phi_n$ for $n \in \NN$, where
$\phi_n(x)=\sqrt{2}\sin(n\pi x)$ and $\lambda_n=-\pi^2 n^2$.

We shall take scalar inputs with $b_n=n \exp(-n^2)$, but, to avoid repeating arguments analogous to those in \cite{jp05,jpp06}, 
we shall consider the case $s=2$, $p>2$,
where Carleson embeddings cannot be tested directly on rectangles,
in order to demonstrate how
Theorem \ref{thm:3.9new} and Remark \ref{rem:312} can be applied.\\

We make use of the following two estimates. 
\begin{enumerate}
\item 
{}From \cite{jp05}, one has
\begin{equation}\label{eq:oldestimate}
1\le  \prod_{j \ne n} \left| \frac{\overline{\lambda_n}+\lambda_j}{\lambda_n-\lambda_j}\right|
\le  \exp(4n(1+\log n))  \qquad \hbox{for each} \quad n \in \NN.
\end{equation}
\item The $L^q(i\RR)$ norm of a Poisson kernel can be estimated by direct integration using (\ref{poisson}),
and one obtains
\begin{equation}\label{eq:normpoisson}
\|p_z\|_{L^q} \asymp (\re z)^{-1+1/q} \qquad \hbox{for} \quad z \in \CC_+,
\end{equation}
which with $q=p/(p-2)$ yields $(\re z)^{-2/p}$.
\end{enumerate}
In the framework of Theorem \ref{thm:3.9new} and Remark \ref{rem:312} with
$\cU=\LL^{-1}H^{p'}(\CC_+,\CC)$ a necessary and sufficient condition for null-controllability
in time $\tau$ is that the function
\begin{equation}\label{eq:24}
\omega \mapsto  \sum_{n \in \NN} 
n^2 \exp(2n^2)\exp(-2n^2\pi^2\tau)
\prod_{j \ne n} \left| \frac{\overline{\lambda_n}+\lambda_j}{\lambda_n-\lambda_j}\right|^2
p_{-\lambda_n}(i\omega)
 \end{equation}
lies in $L^{p/(p-2)}(\RR)$. This is a sufficient condition in the case $\cU=L^p(0,\infty; \CC^N)$, by Corollary~\ref{cor:313}.
Note that checking such a condition is made simpler by the fact that the expression in (\ref{eq:24}) is a sum of positive functions

We deduce easily using (\ref{eq:oldestimate}) and (\ref{eq:normpoisson}) that estimate (\ref{eq:24}) holds if $\tau>1/\pi^2$, since
the series of Poisson kernels converges in $L^{p/(p-2)}$ norm, but does not hold if $\tau<1/\pi^2$, since the series does not converge.
This is in accordance with the results obtained in the case $p=2$.


\section{Conclusions}
\label{sec:4}

We have seen that problems of minimal-norm tangential interpolation can be linked to questions involving Carleson measures and
to more general versions such as those presented in \cite{duren,Lue91}. These in turn have applications to controllability questions
where the input spaces are vectorial  Sobolev spaces or $L^p$ spaces. Provided that the sequence of eigenvalues is reasonably
regularly-distributed, it is possible to solve such questions by the techniques presented above.\\

One significant open question remains, namely, to find an exact necessary and sufficient condition for interpolation in the right half-plane by functions
that are Laplace transforms of $L^p(0,\infty)$ functions; even the discrete case of interpolation  in the disc by an analytic function
whose Fourier coefficients form an $\ell^p$ sequence is only fully solved in the case $p=2$. A full answer to this question would have
immediate applications.

\subsection*{Acknowledgement} The authors gratefully acknowledge support from the Royal Society's {\em International Joint Project\/} scheme.


\end{document}